\documentclass{amsart}
\usepackage{graphicx, amssymb}
\swapnumbers
\newtheorem*{theorem*}{Theorem}
\newtheorem{theorem}{Theorem}[section]
\newtheorem{lemma}[theorem]{Lemma}

\newtheorem{corollary}[theorem]{Corollary}

\theoremstyle{definition}
\newtheorem{definition}[theorem]{Definition}
\newtheorem{example}[theorem]{Example}

\theoremstyle{remark}
\newtheorem{remark}[theorem]{Remark}

\numberwithin{equation}{section}
\newcommand{\firef}[1]{Figure~{\rm\ref{#1}}}
\newcommand{\thref}[1]{Theorem~{\rm\ref{#1}}}

\newcommand{\leref}[1]{Lemma~{\rm\ref{#1}}}
\newcommand{\coref}[1]{Corollary~{\rm\ref{#1}}}

\newcommand{\exref}[1]{Example~{\rm\ref{#1}}}

\newcommand{\seref}[1]{Section~{\rm\ref{#1}}}

\newcommand{\fig}[1]
{\raisebox{-0.5\height}%
{\includegraphics{#1}}}
\newcommand{\st}{\; | \;}                               
\newcommand{\ttt}{\otimes}                               
\newcommand{\tta}{\otimes_A}                               
\newcommand{\tbox}{\boxtimes}

\newcommand{\surjto}{\twoheadrightarrow}      
\newcommand{\injto}{\hookrightarrow}          
\newcommand{\isoto}{\xrightarrow{\sim}}       
\newcommand{\xxto}{\xrightarrow}              

\newcommand{\one}{\mathbf{1}}
\renewcommand{\i}{{\mathrm{i}}}   
\newcommand{\Cset}{\mathbb{C}}    
\newcommand{\Z}{\mathbb{Z}}       
\newcommand{\Rset}{\mathbb{R}}    
\newcommand{\dd}{\mathbf{d}}      
\newcommand{\V}{{\mathcal{V}}}    
\newcommand{\C}{\mathcal{C}}      
\newcommand{\A}{\mathcal{A}}      

\newcommand{\al}{\alpha}
\newcommand{\be}{\beta}

\newcommand{\Ga}{\Gamma}
\newcommand{\de}{\delta}

\newcommand{\ph}{\varphi}
\newcommand{\Ph}{\Phi}

\newcommand{\eps}{\varepsilon}
\renewcommand{\th}{\theta}
\newcommand{\g}{\mathfrak{g}}
\newcommand{\ghat}{\widehat{\mathfrak{g}}}
\newcommand{\U}{U_q(\mathfrak{sl}_2)} 
\newcommand{\slt}{\mathfrak{sl}_2}    
\newcommand{\slthat}{\widehat{\mathfrak{sl}}_2}    

\DeclareMathOperator{\Res}{Res}
\DeclareMathOperator{\Rep}{Rep}
\DeclareMathOperator{\Ind}{Ind}
\DeclareMathOperator{\id}{id}

\DeclareMathOperator{\Hom}{Hom}
\DeclareMathOperator{\Ext}{Ext}

\DeclareMathOperator{\im}{Im}

\begin{document}

\title[On a $q$-analog of McKay correspondence] 
      {On a $q$-analog of the McKay correspondence and 
       the ADE classification of $\slthat$ conformal field theories}

\author{Alexander Kirillov, Jr.}
\address{Department of Mathematics, SUNY at Stony Brook, 
         Stony Brook, NY 11794, USA}
\email{kirillov@math.sunysb.edu}
\urladdr{http://www.math.sunysb.edu/\textasciitilde kirillov/}
\author{Viktor Ostrik}
\address{Department of Mathematics, MIT, 
         Cambridge, MA 02139}
\email{ostrik@math.mit.edu}
\thanks{The first author was supported in part by NSF Grant
  DMS--9970473.\\
    The second author was supported in part by NSF Grant DMS--0098830}
                           
\begin{abstract}
  The goal of this paper is to give a category theory based definition
  and classification of ``finite subgroups in $\U$'' where
  $q=e^{\pi\i/l}$ is a root of unity. We propose a definition of such a
  subgroup in terms of the category of representations of $\U$; we
  show that this definition is a natural generalization of the notion
  of a subgroup in a reductive group, and that it is also related with
  extensions of the chiral (vertex operator) algebra corresponding to
  $\slthat$ at level $k=l-2$.  We show that ``finite subgroups in
  $\U$'' are classified by Dynkin diagrams of types $A_n, D_{2n}, E_6,
  E_8$ with Coxeter number equal to $l$, give a description of this
  correspondence similar to the classical McKay correspondence, and
  discuss relation with modular invariants in $(\slthat)_k$ conformal
  field theory.
  
  The results we get are parallel to those known in the theory of von
  Neumann subfactors, but our proofs are independent of
  this theory.

\end{abstract}

\maketitle
\section*{Introduction}
The goal of this paper is to describe a $q$-analogue of the McKay
correspondence. Recall that the usual McKay correspondence is a
bijection between finite subgroups $\Ga\subset SU(2)$ and affine
simply-laced Dynkin diagrams (i.e., affine ADE diagrams). Under this
correspondence, the vertices of Dynkin diagram correspond to
irreducible representations of $\Ga$ and the matrix of tensor product
with $\Cset^2$ is $2-A_{ij}$ where $A$ is the Cartan matrix of the ADE
diagram (see \cite{mckay}). 

There have been numerous results regarding generalization of McKay
correspondence with $SU(2)$ replaced by $\U$ with $q$ being a root of
unity, $q=e^{\pi\i/l}$.  Of course, since $\U$ is not a group, one
must first find a reasonable way of making sense of this question.
This is usually done by reformulating the problem in terms of
representation theory of $\U$. More precisely, it is known that the
category of finite-dimensional representations of $\U$ has a
semisimple subquotient category $\C$, with simple objects $V_0\dots,
V_k, k=l-2$ (definition of this category was suggested by Andersen and
his collaborators, see \cite{AP}; a review can be found, e.g., in
\cite{Finkelberg} or in \cite{BK}). As was shown by Finkelberg
\cite{Finkelberg}, using results of Kazhdan and Lusztig \cite{KL}`,
the category $\C$ is equivalent to the category of integrable
$\slthat$-modules of level $k$ with fusion tensor product. This latter
category plays a key role in the Wess--Zumino--Witten model of
conformal field theory.  Finally, as was shown by Wassermann and his
students, fusion in this category can also be described in the
language of von Neumann algebras (see, e.g.  \cite{was2},
\cite{laredo}).

Below is an overview of some of the known classification results
related to $\U$, or, more precisely, to category $\C$.

\begin{enumerate}
  
\item {\em Ocneanu's classification of subfactors}.\footnote{We would
    like to thank the referee and M.~M\"uger for explaining to us the
    status of this classification.} It is known that to every
  inclusion of von Neumann factors $N\subset M$ of finite index one
  can associate a number of algebraic structures (index, principal
  graph, relative commutants, etc). In particular, such an inclusion
  defines a tensor category of $N-N$ bimodules. Ocneanu has suggested
  that a subfactor $N\subset M$ with the category of $N-N$ bimodules
  equivalent to $\C$ should be considered as a ``subgroup of $\U$'';
  thus, classification of subgroups reduces to classification of
  subfactors. He also gave \cite{ocneanu} a complete classification of
  such subfactors: they are classified by Dynkin diagrams of types $A$
  (which corresponds to trivial inclusion), $D_{2n}$, $E_6$, $E_8$,
  with Coxeter number equal to $l$. Full proof of this result has been
  given in the works of Popa \cite{popa}, Bion--Nadal \cite{BN},
  \cite{BN2}, Izumi \cite{I}, \cite{I2}.

\item {\em Cappelli-Itzykson-Zuber's classification of modular
    invariants} of conformal field theories based on integrable
  representations of $\slthat$ at level $k=l-2$ (see \cite{CIZ} or the
  review in \cite{CFT}). These modular invariants are classified by
  Dynkin diagrams of ADE type with Coxeter number equal to $l$. It is
  known, however, that modular invariants of types $A, D_{even},
  E_6,E_8$ can be obtained from extensions of the corresponding chiral
  (or vertex operator) algebra, while invariants of the type $E_7,
  D_{odd}$ can not be obtained in this way (see \cite{MST}). This
  classification is related to the previous one: it can be shown that
  every subfactor $N\subset M$ gives a modular invariant, see
  \cite{ocneanu2}, \cite{ocneanu3} and papers of
  B\"ockenhauer--Evans--Kawahigashi \cite{bek, bek2}.

\item {\em Etingof and Khovanov's classification of the ``integer''
  modules over the Gro\-then\-dieck ring} (``fusion algebra'') of $\C$ (see
  \cite{EK}). In this classification, all finite Dynkin diagrams and
  even diagrams with loops appear. 

\end{enumerate}

It should be noted that the classification of ``subgroups in $\U$''
given in the theory of subfactors requires good knowledge of von
Neumann algebras and subfactors. It is very different from the ideas
in the proof of the classical McKay correspondence.

The main goal of the present paper is to study an alternative
definition of a subgroup in $\U$ which uses nothuing but the
theory of tensor categories (which one has to use anyway to work in
$\C$). Namely, a subgroup in $\U$ is by definition a commutative
associative algebra in $\C$, i.e. an object $A\in \C$ with
multiplication morphism $\mu\colon A\ttt A\to A$ satisfying suitably
formulated commutativity, associativity and unit axioms and some mild
technical restrictions.\footnote{While we arrived at this definition
  independently, we are hardly the first to introduce it. This
  definition had also been suggested by Wassermann \cite{was} and
  M\"uger (unpublished).} We argue that this is the right definition for
the following reasons:

\begin{enumerate}
\item If we replace $\C$ by a category of representations of a
  reductive group $G$, then commutative associative algebras in $\C$
  correspond to subgroups of finite index in $G$.
  
\item If we replace $\C$ by a category of representations of some
  vertex operator algebra $\V$ (which is good enough so that $\C$ is a
  modular tensor category, as it happens for all VOA's appearing in
  conformal field theory), then associative commutative algebras in
  $\C$ (with some minor restrictions) exactly correspond to
  ``extensions'' $\V_e\supset \V$ of this VOA; in other words, in this
  way we recover the notion of extension of a conformal field theory.

\item Every subfactor $N\subset M$ defines such an algebra in the
  category of $N-N$ bimodules. 

\end{enumerate}

We show that for any modular category $\C$ a commutative associative
algebra $A\in\C$ gives rise to two different categories of modules
over $A$. One of these categories, $\Rep A$, comes with two natural
functors $F\colon \C\to \Rep A, G\colon \Rep A\to \C$; $F$ is a tensor
functor, so it defines on $\Rep A$ a structure of a module category
over $\C$. There is also a smaller category $\Rep^0 A$; if $A$ is
``rigid'', then both $\Rep A$ and $\Rep^0 A$ are semisimple and
$\Rep^0A$ is modular. Both of these categories have appeared in the
physical literature in the language of extensions of chiral algebra:
in particular, $\Rep^0 A$ is the category of modules over the extended
VOA $\V_e$, and  $\Rep A$ is the category of ``twisted''
$\V_e$-modules. These modules appear as  possible boundary
conditions for extended CFT which preserve $\V$ (see \cite{FS},
\cite{PZ} and references therein); sometimes they are also called
``solitonic sectors''.

Applying this general setup to $\C$ being the semisimple part of
category of representations of $\U$, we see that the fusion algebra of
$\Rep A$ is a module over the fusion algebra of $\C$, which gives a
relation with Etingof-Khovanov classification mentioned above. Using
their results, we get the following classification theorem which we
consider to be the $q$-analogue of McKay correspondence:

\begin{theorem*}
  Commutative associative algebras in $\C$ are classified by the
  \textup{(}finite\textup{)} Dynkin diagrams of the types $A_n,
  D_{2n}, E_6, E_8$ with Coxeter number equal to $l$. Under this
  correspondence, the vertices of the Dynkin diagram correspond to
  irreducible representations $X_i\in \Rep A$ and  the matrix of tensor
  product with $F(V_1)$ in this basis is $2-A$, where $A$ is the Cartan
  matrix of the Dynkin diagram and  $V_1$ is the fundamental
  \textup{(}2-dimensional\textup{)} representation of $\U$. 
\end{theorem*}

Since $\Rep^0 A$ is modular, each of these algebras gives a modular
invariant providing a relation with the ADE classification of
Cappelli-Itzykson-Zuber.

The first part of this theorem --- that is, that commutative
associative algebras are classified by Dynkin diagrams---is hardly
new; in the language of extensions of a chiral algebra, it has been
(mostly) known to physicists long ago (see, e.g.  \cite{MST}), and
these extensions have been studied in a number of papers.  However,
the second part of the theorem, which explicitly describes a
correspondence in a manner parallel to the classical McKay
correspondence to the best of our knowledge is new.

The main result of this theorem is parallel to the classification of
finite subgroups in $\U$ as defined in the theory of subfactors.  Not
being experts in this theory, we are not describing the precise
relation here\footnote{We were recently informed by M.~M\"uger that he
  is currently writing a series of papers, which, among other things,
  will give a detailed review of this connection.}. We just note that
our proofs are completely independent and are not based on the
subfactor theory, even though some of the methods we use (most
notably, the use of conformal embeddings) are parallel to those used
in the subfactor theory.

Finally, it should be noted that the problem of finding $\C$ algebras
$A$ is closely related to the problem of finding all module categories
over $\C$ (such module categories also play important role in CFT; in
physical literature, they are usually described by a certain kind of
$6j$-symbols, see \cite{PZ}).  Indeed, for every $\C$-algebra $A$ the
category $\Rep A$ is a module category over $\C$. It is expected that
in the $\U$ case, all module categories over $\C$ are classified by
all ADE Dynkin diagrams with Coxeter number equal to $l$. Theorem
above gives construction of module categories of type $A_m, D_{even},
E_6, E_8$; it is easy to show that a module category of type
$D_{2n+1}$ can be constructed from representations of some associative
but not commutative $\C$-algebra.  We expect the same to hold for the
module category of type $E_7$.

\subsection*{Note} While working on this paper, we were informed by
A.~Wassermann and H.~Wenzl that they have obtained similar results
based on the subfactor theory. In fact, many of the results of this
paper coincide with those  announced by Wassermann in \cite{was}
(except that we do not use unitary structure), even though we arrived
at these results completely independently.

\subsection*{Note} It was recently pointed out to us that some of the
results about algebras in a braided tensor category and modules over
them which we prove in Section 1 have been previously found in
\cite{Pa}.

\subsection*{Acknowledgments} We would like to thank  Pavel Etingof for many
fruitful discussions and R.~Coquereaux, J.~Fuchs, Y.-Z.~Huang,
Y.~Kawahigashi, J.~Lepowsky, J.~McKay, M.~M\"uger, V.~Petkova,
C.~Schweigert, F.~Xu for comments on the first version of this paper.
We would also like to thank the referee for his detailed comments and
numerous suggestions.

\section{Algebras and modules}\label{sbasic}
Throughout the paper, we denote by $\C$ a semisimple abelian category
over $\Cset$ (most of the results are also valid for any base field
$k$ of characteristic zero). We denote by $I$ the set of isomorphism
classes of irreducible objects in $\C$ and fix some choice of
representative $V_i$ for every $i\in I$. We always assume that the
spaces of morphisms are finite-dimensional; since $\Cset$ is
algebraically closed, this implies that $\Hom_\C(V_i, V_i)=\Cset$. We
will denote $\langle X,Y\rangle= \dim\Hom_\C(X,Y)$; in particular,
$\langle V_i, X\rangle$ is the multiplicity of $V_i$ in $X$ which
shows that this multiplicity is finite. 

We assume that $\C$ is a rigid balanced braided tensor category (see,
e.g., \cite{BK} for a review of the theory of braided tensor
categories). The commutativity isomorphism will be denoted by
$R_{V,W}$; the associativity isomorphism and other canonical
identifications such as $(V\ttt W)^*\simeq W^*\ttt V^*$ will be
implicit in our formulas. Additionally, we require that $\one$ is a
simple object in $\C$. We will use the symbol $0$ to denote the
corresponding index in $I$ : $V_0=\one$.

We denote by $K(\C)$ the complexified Grothendieck ring (``fusion
algebra'') of the category $\C$; this is a commutative associative
algebra over $\Cset$ with a basis given by classes $[V_i]$ of simple
objects.

\begin{definition}\label{dalg}
  An associative commutative algebra $A$ in $\C$ (or $\C$-algebra for
  short) is an object $A\in \C$ along with morphisms $\mu\colon A\ttt
  A\to A$ and $\iota_A\colon \one \injto A$ such that the following
  conditions hold:
\begin{enumerate}
\item (Associativity) 
Compositions $\mu \circ(\mu \ttt \id), \mu\circ (\id\ttt \mu)\colon A^{\ttt
  3}\to A$ are equal. 

\item  (Commutativity) 
Composition $\mu \circ R_{AA} \colon A\ttt A \to A$ is equal to $\mu$.

\item (Unit) 
Composition $\mu \circ (\iota_A\ttt A)\colon A=\one\ttt A\to A$ is equal to
$\id_A$.

\item (Uniqueness of unit)
$\dim \Hom_{\C}(\one, A)=1$.

\end{enumerate}
\end{definition}

The notion of a $\C$-algebra is not new; it had been used in many
papers (for example, in \cite{Ros}). However, most authors only use
algebras in symmetric tensor categories. Algebras in braided
categories were studied in \cite{B}; however, the discussion there is
limited to the case where algebra $A$ is ``transparent'', i.e.
$R_{VA}R_{AV}=\id$ for every $V\in \C$. This setting is too
restrictive for us.  Most of results in \cite{B} generalize to
non-transparent case easily, others (mainly, the results regarding
distinction between two categories of modules, $\Rep A$ and $\Rep^0
A$) require significant work. 

Commutative algebras in a braided tensor category are also discussed
in \cite{Pa}. Unfortunately, this paper was only brought to our
attention after the first version of the current paper appeared in the
electronic archive. For this reason and for  for reader's
convenience, we give here complete proofs of most of the  results; still, we
would like to point out that most results of this section were first
obtained in \cite{Pa} and \cite{B}.

\begin{figure}[h]
\begin{equation*}
\fig{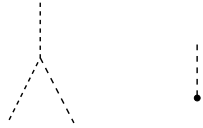}
\end{equation*}
\caption{Morphisms $\mu$ and $\iota_A$}\label{fmuiota}
\end{figure}

We will frequently use graphs to present morphisms in $\C$, as
suggested by Reshetikhin and Turaev. We will use the same conventions
as in \cite{BK}, namely, the morphisms act ``from bottom to top''. We
will use dashed line to represent $A$ and the graphs shown in \firef{fmuiota}
to represent $\mu$ and $\iota_A$. 

\begin{figure}[h]
\begin{equation*}
\fig{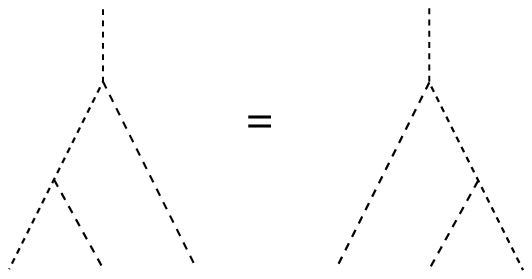}\quad \fig{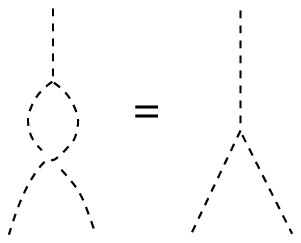}\quad \fig{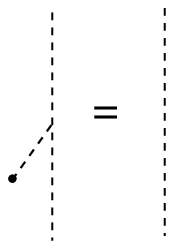}
\end{equation*}
\caption{Axioms of a commutative associative algebra.}\label{fig1}
\end{figure}

With this notation, the axioms of a $\C$-algebra can be presented as
shown in \firef{fig1}. 

We leave it to the reader to define the notions of morphism of
algebras, subalgebras and ideals, quotient algebras etc. 

\begin{figure}[h]
\begin{equation*}
\fig{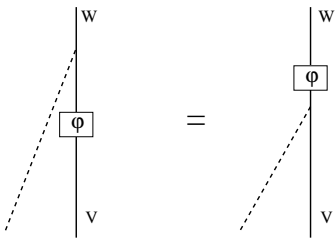}
\end{equation*}
\caption{Definition of morphisms in $\Rep A$}\label{fig2}
\end{figure}

\begin{definition}\label{drepa}
  Let $\C$ be as above and $A$ --- a $\C$ algebra. Define the category
  $\Rep A$ as follows: objects are pairs $(V, \mu_V)$ where $V\in \C$
  and $\mu_V\colon A\ttt V\to V$ is a morphism in $\C$ satisfying the
  following properties: 
\begin{enumerate}
\item
$\mu_V \circ (\mu\ttt \id)=\mu_V \circ (\id \ttt
  \mu_V)\colon A\ttt A\ttt V\to V$
\item $\mu_V(\iota_A \ttt \id)=\id\colon \one\ttt V\to V$
\end{enumerate}

The morphisms are defined by 
\begin{multline}\label{amorph}
\Hom_{\Rep A}((V, \mu_V), (W, \mu_W))\\
=\{ \ph\in \Hom_\C(V, W)\st 
                   \mu_W\circ(\id \ttt \ph)=\ph\circ\mu_V\colon A\ttt V\to
                   W\} 
\end{multline}
(see \firef{fig2}).

\end{definition}

An instructive example of such a situation is when $G$ is a finite
group and $\C$ is the category of finite-dimensional complex
representations of $G$. In this case we will show that semisimple
$\C$-algebras correspond to subgroups in $G$ (see \seref{sgroup}).

\begin{remark}\label{rintuition}
Contrary to the usual intuition, typically the larger $A$, the
smaller is its category of representations. In the above mentioned
example $\C=\Rep G$, correspondence between subgroups $H\subset G$ and
$\C$-algebras is given by $A=F(G/H)$, so large $A$ corresponds to
small $H$ and thus, to small $\Rep A=\Rep H$.
\end{remark}

Let us study basic properties of $\Rep A$. For brevity, we will use
notation $\Hom_A$ instead of $\Hom_{\Rep A}$.

\begin{lemma}\label{labelian}\ \\
\begin{enumerate}
\item $\Rep A$ is an abelian category with finite-dimensional spaces
  of morphisms; every object in $\Rep A$ has finite length. 
\item $\Hom_A(A,A)=\Cset$. 
\end{enumerate}
\end{lemma}
\begin{proof}
Since $\C$ is an abelian category, it suffices to prove that for
$f\in \Hom_{A}(V, W)$, $\im f$ and $\ker f$ are actually
$A$-submodules in $W, V$ respectively. The check is straightforward
and is left to the reader. 

Let $\ph\in \Hom_A(A,A)$. By definition we have:
\begin{equation*}
\ph = \ph \mu (\id_A\ttt \iota_A)=\mu (\id \ttt \ph \iota_A)
\end{equation*}
But since $\one$ has multiplicity one in $A$, one has $\ph \iota_A=c \iota_A$
for some constant $c$. Thus, $\ph = c\mu (\id_A\ttt \iota_A)=c\id$. 
\end{proof}

\begin{theorem}\label{ttensor}
$\Rep A$ is a monoidal category with unit object $A$.

\end{theorem}
\begin{proof}
Let $V, W\in \Rep A$. Define $V\tta W=V\ttt W/\im(\mu_1-\mu_2)$ where
$\mu_1, \mu_2\colon A\ttt V\ttt W\to V\ttt W$ are defined by 
\begin{align*}
\mu_1&=\mu_V\ttt\id_W,\\
\mu_2&=(\id_V\ttt \mu_W)R_{AV}.
\end{align*}
This defines $V\tta W$ as an object of $\C$. Define $\mu_{V\tta W}$ to
be $\mu_1$ or $\mu_2$ which obviously give the same morphism. One
easily sees that this defines a structure of $A$-module on $V\ttt W$
and that so defined tensor product is associative. To
check that $A$ is the unit object, consider morphisms $\mu_V\colon
A\tta V\to V$ and $\iota_A\ttt\id_V\colon  V\to A\tta V$. Straightforward
check shows that they are well defined, commute with the action of $A$
(that is, satisfy \eqref{amorph}) and thus define morphisms in $\Rep
A$ and finally, that they are inverse to each other. 
\end{proof}

\begin{theorem}\label{tfunctors}
Define functors $F\colon \C\to \Rep A, G\colon \Rep A \to \C$ by 
$F(V)=A\ttt V, \mu_{F(V)}=\mu \ttt \id$ and $G(V, \mu_V)=V$. Then 

\begin{enumerate}
\item Both $F$ and $G$ are exact and injective on morphisms.
\item $F$ and $G$ are adjoint: one has canonical functorial
  isomorphisms 
$$
\Hom_A(F(V), X)=\Hom_\C(V, G(X)), \qquad V\in \C, X\in \Rep A.
$$

\item $F$ is a tensor functor: one has canonical isomorphisms $F(V\ttt
  W)= F(V)\tta F(W), F(\one)=A$. 
\item One has canonical isomorphisms $G(F(V))=A\ttt V$ and, more
  generally, $G(F(V)\tta X)=  V\ttt G(X)$. 
\end{enumerate}
\end{theorem}

\begin{proof}
 Part (1) is obvious; for part (2), define maps $\Hom_A(F(V), X)\to
 \Hom_\C(V, G(X))$ and $\Hom_\C(V, G(X))\to \Hom_A(F(V), X)$ as shown
 in \firef{fadjoint1}; it is easy to deduce from the axioms that these
 maps are inverse to each other. 

\begin{figure}[h]
\begin{equation*}
\fig{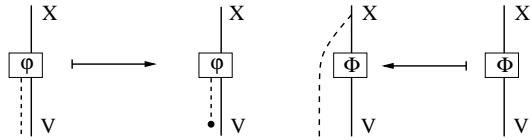}
\end{equation*}
\caption{Identifications $\Hom_A(F(V), X)=\Hom_\C(V, G(X))$.
        Here  $\ph\in\Hom_A(F(V), X), \ \Phi\in\Hom_\C(V, G(X))$.}
        \label{fadjoint1}
\end{figure}

  To prove that $F$ is a tensor functor, define functorial morphisms
  $f\colon F(V\ttt W)\to F(V)\tta F(W), g\colon F(V)\tta F(W)\to
  F(V\ttt W)$ by 
\begin{align*}
f&=\id_A\ttt\id_V\ttt \iota_A\ttt \id_W\colon A\ttt V\ttt
  W\to (A\ttt V)\tta (A\ttt W)\\
g&\colon (A\ttt V)\tta (A\ttt W)\xxto{R_{AV}^{-1}} A\tta
  A\ttt V \ttt V \xxto{\mu} A\ttt V\ttt W.
\end{align*}
It is immediate to check that they are well-defined and inverse to
each other. 
\end{proof}

\begin{corollary}\label{cmodule}
$\Rep A$ is a module category over $\C$, i.e. there is an additive  functor
$\tbox\colon \C\times \Rep A\to \Rep A$ and  isomorphisms 
\begin{align*}
&(V_1\ttt V_2)\tbox X\simeq V_1\tbox (V_2\tbox X), 
            \qquad V_1, V_2\in \C, X\in \Rep A\\
&\one\tbox X\simeq X, \quad X\in \Rep A
\end{align*}
satisfying usual compatibility conditions.
\end{corollary}
\begin{proof}
Suffices to take $V\tbox X=F(V)\tta X$. 
\end{proof}

In particular, this implies that the Grothendieck group $K(\Rep A)$ is
a module over the Grothendieck ring $K(\C)$.

However, it is not true that $\Rep A$ is a braided tensor
category. In order to get a braided structure, we need to consider a
smaller category.

\begin{definition}\label{drep0}
$\Rep^0A$ is the full subcategory in $\Rep A$ consisting of objects
$(V, \mu_V)$ such $\mu_V\circ R_{VA}R_{AV}=\mu_V$. 
\end{definition}

\begin{remark}
In \cite{Pa}, such $A$-modules are called ``dyslectic''. 
\end{remark}

If $\C$ is symmetric, then $\Rep^0 A=\Rep A$. More generally, the same
holds if  $A$ is ``transparent'',
or ``central'', in $\C$ (that is, $R_{VA}R_{AV}=\id$ for every $V\in
\C$); this is the situation considered in \cite{B}. However, in many
interesting cases $A$ is not central, and $\Rep A\ne \Rep^0 A$. 

Later we will justify this definition by showing that if $\C$ is a
category of representation of some vertex operator algebra, and $A$ is
an extended vertex operator algebra, then $\Rep^0A$ (and not $\Rep
A$!) is exactly the category of representations of the vertex operator
algebra $A$.

\begin{theorem}[\cite{Pa}]
\label{tbraided}
The category $\Rep^0 A$ is a braided tensor category, with the
commutativity isomorphism inherited from $\C$.
\end{theorem}
\begin{proof}
  Let us first show that for $X,Y\in \Rep^0 A$, 
  $X\tta Y\in \Rep^0 A$. This follows from the sequence of identities
  shown in \firef{frtta}. The notation $f_1\equiv f_2$ for $f_1, f_2\colon
  A\ttt X\ttt Y\to X\ttt Y$ means that the induced operators $A\ttt
  (X\ttt Y)\to X\tta Y$ are equal, i.e.  $p\circ f_1=p\circ f_2$,
  where $p\colon X\ttt Y\to X\tta Y$ is the canonical projection.  
\begin{figure}[ht]
$$
\fig{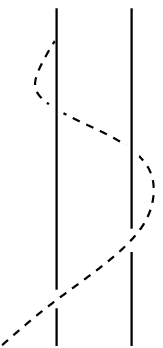}\quad=\quad\fig{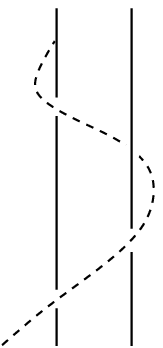}\quad\equiv\quad\fig{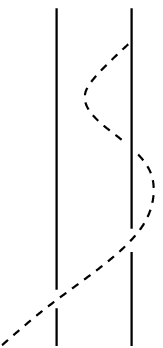}
\quad=\quad\fig{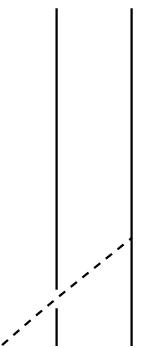}\quad\equiv\quad\fig{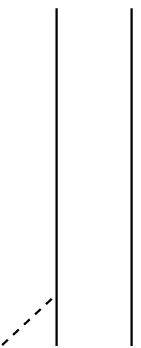}
$$
\caption{}\label{frtta}
  \end{figure}

  Next, we need to show that the commutativity isomorphism
  $R_{XY}\colon X\ttt Y\to Y\ttt X$ descends to isomorphism $ X\tta
  Y\to Y\tta X$.  This is equivalent to showing that $R_{XY}
  (I)\subset I$, where $I=\im (\mu_1-\mu_2)$ is the kernel of the
  canonical projection $X\ttt Y\to X\tta Y$.  To do so, let us rewrite
  composition $R_{XY}\circ (\mu_1-\mu_2)$ as shown in
  \firef{frcommut}.  Thus, $R_{XY}(\mu_1-\mu_2)\equiv 0$, or,
  equivalently, $R_{XY}(I)\subset I$.
\begin{figure}[ht]
\begin{align*}
& \fig{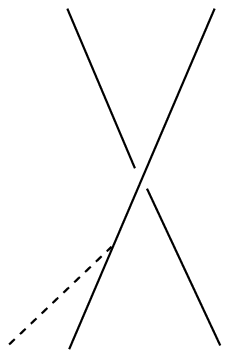}-\fig{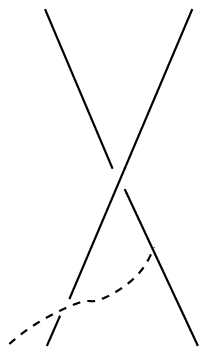}
  =\fig{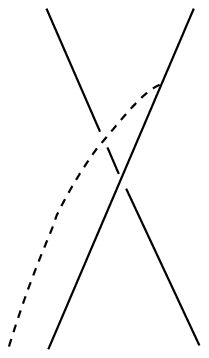}-\fig{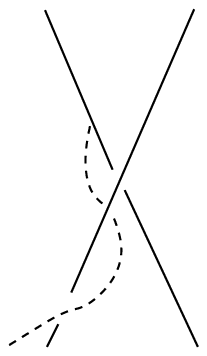}\\
&\qquad
  \equiv \fig{rcommutc}-\fig{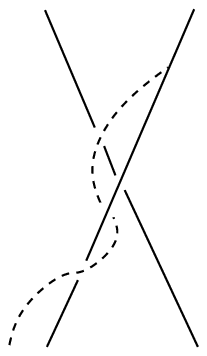}\equiv 0.
\end{align*}
\caption{}\label{frcommut}
  \end{figure}
 
Abusing the language, we will also use the same notation $R_{X,Y}$ for
the descended morphisms $X\tta Y\to Y\tta X$. Then it is immediate
from the definition that these morphisms are $A$-morphisms. Finally,
since the commutativity isomorphism in $\C$ satisfies the hexagon
axioms, the same must hold for the descended operators; thus, the
descended operators $R_{X,Y}\colon X\tta Y\to Y\tta X$ define a
structure of a braided tensor category on $\Rep^0 A$. 
\end{proof}

Analysis of this proof also shows the reason why the larger category
$\Rep A$ is not braided: the last identity in \firef{frcommut} would
fail.

\begin{figure}[h]
\begin{equation*}
\fig{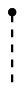}
\end{equation*}
\caption{Morphism $\eps_A$}\label{feps}
\end{figure}

We also need to know whether categories $\Rep A, \Rep^0 A$ are rigid.
Define $\eps_A\colon A\to \one$ so that $\eps_A \iota_A=\id$ (recall
that $\langle A,\one\rangle=1$, so this condition uniquely defines
$\eps_A$). We will use the graph shown in \firef{feps} to represent
$\eps_A$.

We will say that a $\C$-morphism $f\colon V\ttt W\to \one$ defines a
{\em non-degenerate pairing} if $f$ defines an isomorphism $V\simeq
W^*$, i.e. there exists a map $g\colon \one\to W\ttt V$ such that
$f,g$ satisfy the rigidity axioms.

\begin{definition} \label{drigid}
A $\C$-algebra $A$ is called {\em rigid} if the map 
\begin{equation}\label{erigid}
e_A\colon A\ttt A\xxto{\mu}A\xxto{\eps_A}\one
\end{equation}
is a non-degenerate pairing  and $\dim_\C A\ne 0$. 
\end{definition}

If $A$ is rigid, then there is a unique morphism $i_A\colon \one\to
A\ttt A$ such that $e_A, i_A$ satisfy the rigidity axioms: the
compositions $A\xxto{\id\ttt i_A}A\ttt A\ttt A\xxto{e_A\ttt \id}A$,
$A\xxto{i_A\ttt \id }A\ttt A\ttt A\xxto{\id \ttt e_A}A$ are equal to
identity (uniqueness follows from a well-known fact that the dual
object is unique up to a unique isomorphism).  We will frequently use
the following simple lemma.

\begin{lemma}\label{lsymmetry}
Both $e_A$, $i_A$ are symmetric: 
\begin{align*}
&e_A R_{AA}=e_A,\\
&R_{AA}i_A=i_A.
\end{align*}
\end{lemma}
\begin{proof} 
  The identity for $e_A$ immediately follows from commutativity of
  multiplication. For $i_A$, it suffices to prove that $R_{AA}i_A$
  satisfies the rigidity axioms, i.e.  both compositions below are
  equal to identity
\begin{align*}
&A\xxto{\id\ttt R_{AA}i_A}A\ttt A\ttt A\xxto{e_A\ttt \id}A,\\
&A\xxto{R_{AA}i_A\ttt \id }A\ttt A\ttt A\xxto{\id \ttt e_A}A.
\end{align*}
To prove the first identity, we represent the
composition by a graph and manipulate it as shown in \firef{ficommut}.
The second identity is proved in the same manner.

\begin{figure}[h]
\begin{equation*}
\fig{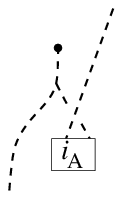}=
\fig{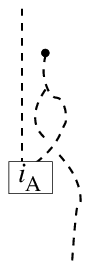}=
\fig{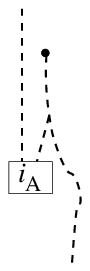}=\id_A
\end{equation*}
\caption{Proof of $R_{AA}i_A=i_A$.}\label{ficommut}
\end{figure}
\end{proof}

Finally, we also need to discuss the following subtle point. 
The map $e_A$ allows us to identify $A\simeq A^*$ and thus, we can
also identify $A\simeq A^{**}$. On the other hand, in any balanced
rigid braided category one has a canonical identification $\de_V\colon
V\simeq V^{**}$. It is natural to ask whether these two
identifications coincide. The answer is given by the following lemma.

\begin{lemma}\label{tha=1}
  Let $A$ be a rigid algebra. Then the morphism $A\simeq A^{**}$
  defined by $e_A$ coincides with $\de_A\colon A\to A^{**}$ iff
  $\th_A=\id$. 
\end{lemma} 
\begin{proof}
Recalling the relation between the twists $\th_V$ and $\de_V$ (see,
e.g., \cite[Section 2.2]{BK}), we see
that the statement of the lemma is equivalent to the following
equation:
$$
\fig{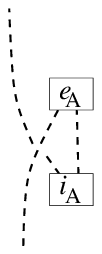}=\th_A.
$$
But it easily follows from symmetry of $e_A$ that the right hand side
is the identity morphism. 
\end{proof}
This lemma shows that if $A$ is rigid, $\th_A=\id$, then  we can
identify $A\simeq A^*$ so that the canonical  morphisms $A\ttt
A^*\to \one, A^*\ttt A\to \one$ are both given by  $e_A$, and the
morphisms $\one\to A\ttt A^*, \one \to A^*\ttt A$ are both given by
$i_A$. As usual, we will use ``cap'' and ``cup'' to denote $e_A,
i_A$  in the figures. Then the statement of \leref{lsymmetry} can be
graphically presented  as follows: 

\begin{equation}\label{esymmetry}
\fig{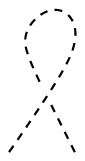}=\fig{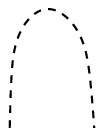}\qquad\qquad
\fig{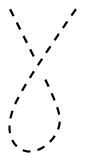}=\fig{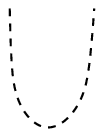}
\end{equation}

We will frequently use the following easy lemma.
\begin{lemma}\label{lnoose}
If $A$ is a rigid $\C$-algebra, $\th_A=\id$, then 
\begin{equation}
\fig{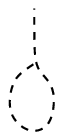}=\dim A \fig{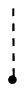}
\end{equation}
\end{lemma}
The proof is immediate if we note that both sides are morphisms $\one
\to A$ and by uniqueness of unit axiom must be proportional.

\begin{theorem}\label{trigid}
  If $\C$ is a rigid category, $A$ --- a rigid $\C$-algebra,
  $\th_A=\id$, then the categories $\Rep A, \Rep^0 A$ are rigid.
\end{theorem}
\begin{proof}
Let $(V, \mu_V)\in \Rep A$. Define the dual object $(V^*, \mu_{V^*})$
as follows: $V^*$ is the dual of $V$ in $\C$ and $\mu_{V^*}$ is defined
by \firef{fvdual}.

\begin{figure}[h]
\begin{equation*}
\fig{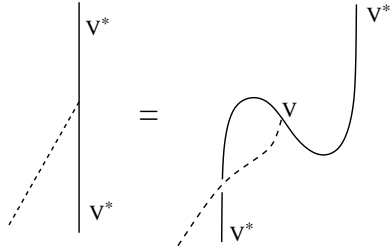}
\end{equation*}
\caption{Definition of dual object in $\Rep A$.}\label{fvdual}
\end{figure}

This definition implies the following identities:

\begin{figure}[h]
\begin{equation*}
\fig{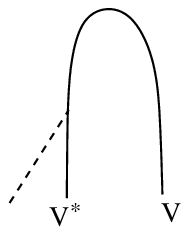}=\fig{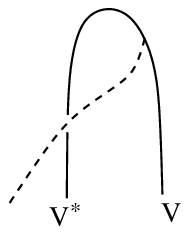}
\qquad \fig{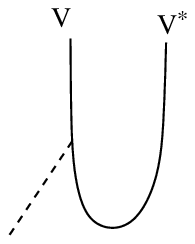}=\fig{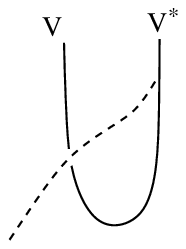} 
\end{equation*}
\caption{}\label{fdualident}
\end{figure}

Define now the maps $\tilde i_V\in\Hom_{A}(A, V\tta V^*), \tilde e_V\in
\Hom_A(V^*\tta  V, A)$ by \firef{frigid} (we leave it to the reader to
check that these formulas indeed define morphisms in $\Rep A$).

\begin{figure}[h]
\begin{equation*}
\tilde i_V=\frac{1}{\dim A}\fig{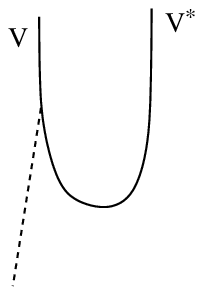}\qquad
\tilde e_V=\fig{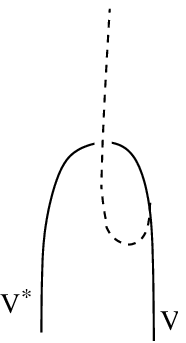}
\end{equation*}
\caption{Rigidity maps in $\Rep A$}\label{frigid}
\end{figure}

It is easy to check by  using identities in \firef{fdualident} 
and isomorphisms $A\tta V\simeq V$ defined in the proof of \thref{ttensor}
that these two maps satisfy the rigidity axioms. 
\end{proof}

\begin{lemma}
Let $A$ be a rigid $\C$-algebra. Then
\begin{enumerate}
\item $F$ and $G$ are 2-sided adjoints of each other: in addition to
  results of \thref{tfunctors}, we also have canonical isomorphisms  
$$
\Hom_A(X, F(V))=\Hom_\C(G(X),V), \qquad V\in \C, X\in \Rep A.
$$
\item 
$F(V^*)\simeq (F(V))^*$.
\end{enumerate}
\end{lemma}

\begin{proof}
To prove part (1), we construct linear maps between $\Hom_A(X, F(V))$
and $\Hom_\C(G(X),V)$ as shown in \firef{fadjoint2}; we leave it to
the reader to check that these maps are inverse to each other.

\begin{figure}[h]
\begin{equation*}
\fig{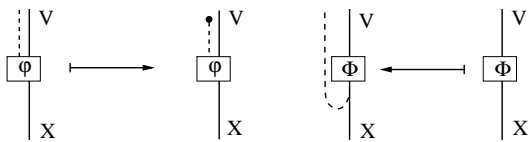}
\end{equation*}
\caption{Identifications $\Hom_A(X, F(V))=\Hom_\C(G(X), V)$.
        Here  $\ph\in\Hom_A(X, F(V)), \ \Phi\in\Hom_\C(G(X), V)$.}
        \label{fadjoint2}
\end{figure}

To prove (2), note that as object of $\C$, $(F(V))^*=V^*\ttt
A^*=V^*\ttt A$, where we used rigidity to identify $A=A^*$. Consider
the morphism $R_{AV^*}\colon A\ttt V^*\to V^*\ttt A$. Again, we leave
it to the reader to check that this morphism is actually a morphism of
$A$-modules $F(V^*)\to (F(V))^*$. 

This shows that $\Rep A$ is rigid. To prove rigidity of $\Rep^0 A$, it
suffices to show that for $X\in \Rep^0 A, X^*\in \Rep^0 A$ which
easily follows from the definition of $\Rep^0 A$ and the definition of
dual object given by \firef{fvdual}. 
\end{proof}

Finally, we need to check discuss whether $\Rep A, \Rep^0 A$ are
balanced. Recall that 
balancing in a rigid braided tensor category is a system of functorial
isomorphisms $\de_V\colon V\simeq V^{**}$ satisfying conditions 
\begin{equation}\label{edev}
\begin{aligned}
{}&\de_{V\ttt W}=\de_V\ttt \de_W,\\
&\de_\one =\id,\\
&\de_{V^*}=(\de_V^*)^{-1},
\end{aligned}
\end{equation}
where we have used canonical identifications $(V\ttt W)^*=W^*\ttt V^*$
and for $F\colon X\to Y$, $f^*\colon Y^*\to X^*$ is the adjoint
morphism. 

This is equivalent to defining a system of functorial morphisms $\th_V
\colon V\to V$ (twists), satisfying 
\begin{equation}\label{ebalancing}
\begin{aligned}
{}&\th_{V\ttt W}=R_{WV}R_{VW}\th_V\ttt \th_W,\\
&\th_\one=\id,\\
&\th_{V^*}=(\th_V)^*.
\end{aligned}
\end{equation}
(see, for example,  \cite[Section 2.2]{BK}).

\begin{theorem}\label{tbalanced}
  Let $\C$ be a rigid balanced braided category, and $A$---a rigid
  $\C$-algebra, $\th_A=\id$. Then

\begin{enumerate}
\item $\Rep^0 A=\{V\in \Rep A\st \th_V \text{ is an $A$-morphism }\}$.
\item $\Rep^0A$ is a rigid balanced braided category, with $\th$
  inherited from $\C$.
\item For any $V\in \Rep A$, the morphism $\de_V\colon V\to V^{**}$
  is an $A$-morphism.  
\end{enumerate}
\end{theorem}
\begin{proof}
  Part (1) follows from $\th_{A\ttt V}=R_{VA}R_{AV}\th_A\ttt \th_V$
  and $\th_A=\id$; (2) immediately follows from (1). 
  To prove (3), it suffices to prove that $\de_V^{-1}\colon V^{**}\to V$ is an
  $A$-morphism. We can rewrite $\de_V^{-1}$ in terms of $\th$ as
  follows (see \cite[Section 2.2]{BK}): 
$$
\de_V^{-1}=\fig{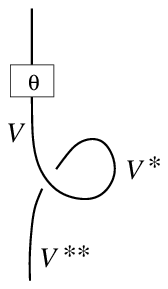}
$$

It now follows from the identities shown in \firef{fdev} 
 (which uses \eqref{ebalancing} and
identities from \firef{fdualident}) that $\de^{-1}$ is an
$A$-morphism. 

\begin{figure}[ht]
$$
\fig{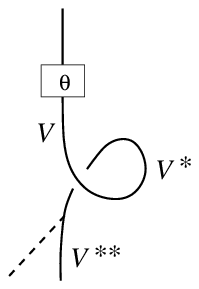}=\fig{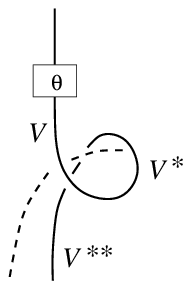}=\fig{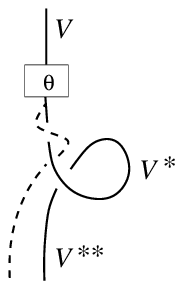}=\fig{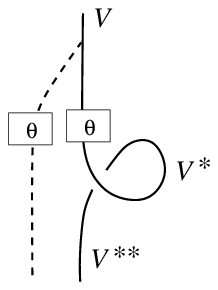}
$$
\caption{Proof that $\de_V^{-1}$ is an $A$-morphism}\label{fdev}
\end{figure}

\end{proof}

This theorem shows that the category $\Rep A$, which is a rigid
monoidal category, while not braided, does have a system of functorial
morphisms $\de_X\colon X\to X^{**}$ satisfying \eqref{edev}. Such
categories are sometimes called ``pivotal''; in such a category, one
can define for every object two numbers, its ``left'' and ``right''
dimension   (see, e.g., \cite{BW}).  We will denote by 
$\dim_A X$ the ``left'' dimension of an object $X\in \Rep A$.

\begin{theorem}\label{tdim0}
  Let $\C$ be a rigid balanced braided category, and $A$---a rigid
  $\C$-algebra such that $\th_A=\id_A$. Then for every $X,Y \in \Rep
  A$, $\dim_{A}(X\tta Y)=\dim_{A}(X)\dim_{A}(Y)$ and 
\begin{align*}
&\dim_{A}(X)=\frac{\dim_{\C}(X)}{\dim_{\C}A},\\
&\dim_{A}(F(V))=\dim_\C(V).
\end{align*}
\end{theorem}
\begin{proof}
Formula $\dim_{A}(X\tta Y)=\dim_{A}(X)\dim_{A}(Y)$ holds in any
pivotal category and can be easily deduced from \eqref{edev}.

Using definition of rigidity morphisms in $\Rep A$ shown in
\firef{frigid},  we see that
$\dim_A X$ is defined by the following identity:
$$
\frac{1}{\dim A}\fig{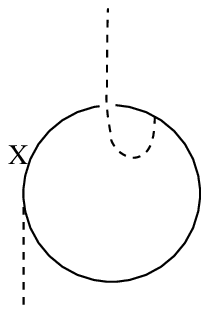} =(\dim_A X) \id_A.
$$
Both sides are $A$-morphisms $A\to A$. Composing them with $\iota_A,
\eps_A$,  we get $\frac{1}{\dim A}\dim_\C X=\dim_A X$. 

Applying this to $X=F(V)=V\ttt A$, we get 
$$
\dim_A F(V)=\frac{1}{\dim
  A} \dim_\C F(V)= \frac{1}{\dim A} (\dim  A)\cdot(\dim V)=\dim V.
$$ 
\end{proof}

As a useful corollary, we get the following result:
\begin{equation}\label{dimA}
\dim_\C(X\tta Y)=\frac{\dim_\C(X)\dim_\C(Y)}{\dim_\C(A)}.
\end{equation}

\begin{remark}
In the theorem above, we could have used ``right'' dimension instead of
the ``left'' dimension (this would require minor change in the
proof). Thus, we see that both left and right dimension of $X\in \Rep
A$ are equal to  $\frac{\dim_{\C}(X)}{\dim_{\C}A}$; in particular,
they are equal to each other. In a similar way, one could prove that
for each $A$-morphism $F\colon X\to X$, its left and right dimension
are equal; in terminology of \cite{BW}, $\Rep A$ is a {\em spherical}
category. 
\end{remark}

For future use, we note the following somewhat unusual result. 
\begin{lemma}\label{lideals}
  Let $\C$ be rigid and $A$ -- a $\C$-algebra such that $\th_A=\id,
  \dim_\C A\ne 0$. Then $A$ is a rigid $\C$-algebra iff $A$ is simple
  as an $A$-module.
\end{lemma}
\begin{proof}
  Let $A$ be rigid; assume $I\subset A$ is a submodule. By rigidity,
  $\one\subset \mu(A\ttt I)$. On the other hand, since $I$ is a
  submodule, this implies that $\one \subset I$. By unit axiom, this
  implies $I=A$.

Conversely, assume that $A$ is simple as $A$-module. Consider $A^*\in\C$
and define on it the action of $A$ as in \thref{trigid}. Then one
easily sees that the morphism 
$$
A\xxto{\id\ttt i_A} A\ttt A\ttt A^*\xxto{e_A\ttt \id}A^*,
$$
where $e_A$ is as in \eqref{erigid} is a morphism of $A$-modules.
On the other hand, usual arguments show that if $A$ is a simple
$A$-module, then so is $A^*$. Thus, such a map is either zero
(impossible because of the  unit axiom) or an isomorphism. 
\end{proof}

\begin{figure}[h]
\begin{equation*}
Q=\frac{1}{\dim A}\quad\fig{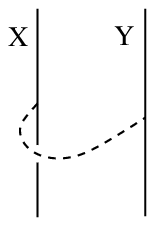}
\end{equation*}
\caption{Projector on $X\tta Y\subset X\ttt Y$}\label{fQ}
\end{figure}

Finally, recall that we defined $X\tta Y$ as a quotient of $X\ttt
Y$. It turns out that in the rigid case, $X\tta Y$ can also be
described as a sub-object of $X\ttt Y$ and thus, as a direct summand. 
\begin{lemma}\label{lsummand}
  If $A$ is a rigid $\C$-algebra, $X,Y\in \Rep A$ and $Q\colon X\ttt
  Y\to X\ttt Y$ is as shown in \firef{fQ}, then $Q^2=Q$ and $\ker
  Q=\ker(X\ttt Y\to X\tta Y)$.
\end{lemma}
Proof of this lemma is left to the reader as an exercise.

\begin{corollary}
  If $A$ is a rigid $\C$-algebra, then one has a canonical direct sum
  decomposition $X\ttt Y=Z\oplus X\tta Y$ for some $Z\in \C$.
\end{corollary}
Indeed, it suffices to take $Z=\ker Q$ and identify $X\tta Y=\im Q$.

\section{Example: groups and subgroups}\label{sgroup}
In this section we discuss an important example of the general setup
discussed above. Namely, let $G$ be a group such that the category
$\C$ of finite-dimensional complex representations of $G$ is
semisimple (for example, $G$ is a finite group or $G$ is a reductive
Lie group). Also, for a finite set $X$ let $F(X)$ be the space of
complex-valued functions on $X$.

\begin{theorem} 
  If $H\subset G$ is a subgroup of finite index, then
  the space $A=F(G/H)$ of functions on $G/H$ is a semisimple
  $\C$-algebra and $\Rep A$ is equivalent to the category $\Rep H$ of
  representations of $H$; under this equivalence the functors $F$ and
  $G$ are identified with the restriction and induction functor
  respectively:
\begin{align*}
F=\Res^G_H\colon& \Rep G \to \Rep H,\\
G=\Ind_H^G\colon& \Rep H\to \Rep G
\end{align*}

\end{theorem}
\begin{proof}
  By definition, an object of $\Rep A$ is a $G$-module $V$ with a
  decomposition $V=\oplus_{x\in G/H} V_x$ such that $gV_x=V_{gx}$, and
  tensor product in $\Rep A$ is given by $(V\tta W)_x=V_x\ttt W_x$.
  Define functor $\Rep A\to \Rep H$ by $\oplus V_x\mapsto V_1$ and
  $\Rep H\to \Rep A$ by $E\mapsto \Ind_H^G E$ (note that it follows
  from definition of the induced module that $V=\Ind E$ has a natural
  decomposition $V=\oplus_{x\in G/H}V_x$). It is trivial to check that
  these functors preserve tensor product and  are inverse to each
  other. 
\end{proof}

\begin{theorem}
For $\C=\Rep G$, any rigid $\C$-algebra is of the form $F(G/H)$
for some subgroup $G$ of finite index. 
\end{theorem}
\begin{proof}
  First, a $\C$-algebra is just a commutative associative algebra over
  $\Cset$ on which $G$ acts by automorphisms. Next, if $A$ is rigid,
  then $A$ is semisimple as a commutative associative algebra over
  $\Cset$.  Indeed, let $N$ be the radical of $A$; then $N$ is
  invariant under the action of $G$ and thus is an ideal in $A$ in the
  sense of $\C$-algebras. By \leref{lideals}, $N=0$. 
 
  Thus, $A$ is the algebra of functions on a finite set $X$ (which can
  be described as the set of primitive idempotents of $A$) and $G$
  acts by permutations on $X$. Since $\Cset$ appears in decomposition
  of $A$ as $G$-module with multiplicity one, this implies that the
  action of $G$ on $X$ is transitive, so $X=G/H$.
\end{proof}

\section{Semisimplicity}
As before, we let $\C$ be a braided tensor category. 

\begin{definition}\label{dsemisimple}
A $\C$-algebra is called semisimple if $\Rep A$ is semisimple.
\end{definition}


We will be mostly interested in the case when $\C$ is rigid and
balanced. In this case, semisimplicity of  $\Rep A$ implies
semisimplicity of  $\Rep^0 A$. 

\begin{theorem}\label{tssrep0}
Let $\C$ be rigid balanced, and $A$ --- a semisimple
$\C$-algebra with $\th_A=\id$. Then 
\begin{enumerate}
\item If $X\in \Rep A$ is simple, then $X\in \Rep^0 A$ iff
  $\th_X=c\cdot \id$. 

\item $\Rep^0 A$ is semisimple, with simple objects $X_{\pi}$ where
  $X_\pi$ is a simple object in $\Rep A$ such that $\th_X=c\cdot \id$.
\end{enumerate}
\end{theorem}
\begin{proof}
Immediately follows from \thref{tbalanced} and the fact that for a
simple object $X\in \Rep A$, $\Hom_A(X,X)=\Cset$. 
\end{proof}

The main result of this section is the following theorem.

\begin{theorem}\label{trigidsemisimple}
Let $\C$ be rigid, and $A$ --- a rigid $C$-algebra. Then $A$  is semisimple. 
\end{theorem}

\begin{proof}
The proof is based on the following lemma.
\begin{lemma}[\cite{B}]
If $A$ is rigid, then every $X\in \Rep A$ is a direct summand in
$F(V)$ for some $V\in \C$.
\end{lemma}
\begin{proof}[Proof of the lemma]
Consider the map $\mu\colon A\ttt A\to A$. It is surjective and
is a morphism of $A$-modules. Moreover, both $A$ and $A\ttt A$ have
canonical structures of $A$-bimodules, and $\mu$ is a morphism of
$A$-bimodules (we leave the definition of $A$-bimodule as an exercise
to the reader). This map has one-sided inverse: if we define
$\ph\colon A\to A\ttt A$ by
$$
\ph=\frac{1}{\dim A}\fig{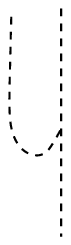}
$$
then $\ph$ is a morphism of $A$-bimodules and it immediately
follows from \leref{lnoose} that $\mu\ph=\id_A$. Thus, $A\ttt A$
splits: we can write 
$$
A\ttt A\simeq A\oplus Z
$$ 
for some $A$-bimodule $Z$ so that under this isomorphism, $\mu$ is the
projection on the first summand. 

Therefore, $A\tta X\simeq X$ is a direct summand of 
$(A\ttt A)\tta X= A\ttt (A\tta X)=A\ttt X=F(G(X))$. 
\end{proof}

From this lemma, the proof is easy. Indeed, it easily follows from
exactness of $G$ and adjointness of $F$ and $G$ (see
\thref{tfunctors}) that for every $V\in \C$, $F(V)$ is a projective
object in $\Rep A$. Since a direct summand of a projective object  is
projective, the lemma implies that every $X\in \Rep A$ is projective
and thus, $\Ext^1(X,Y)=0$ for every $X, Y\in \Rep A$.
\end{proof}

\begin{remark}
  Morally, this theorem is parallel to the following well known result
  in Lie algebra theory: if the Killing form on $\g$ is non-degenerate,
  then the category of finite-dimensional representations of $\g$ is
  semisimple (this is combination of Cartan's criterion of
  semisimplicity and Weyl's complete reducibility theorem). The proof,
  of course, is completely different.
\end{remark}

\section{Modularity}
\begin{figure}[h]
\begin{equation*}
\fig{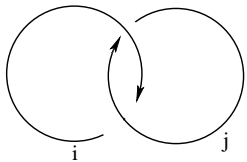}
\end{equation*}
\caption{Matrix $\tilde s_{ij}$}\label{fsij}
\end{figure}

Recall that a semisimple balanced rigid braided category $\C$ is
called modular if it has finitely many isomorphism classes of simple
objects $V_i, i\in I, |I|<\infty$  and the matrix $\tilde s_{ij}$
defined by \firef{fsij} is non-degenerate. (We will also use numbers
$\tilde s_{VW}$ defined in the same way as $\tilde s_{ij}$ but with
$V_i$ replaced by $V$, $V_j$ replaced by $W$).

In this case, it is known that the matrices 
\begin{equation}\label{st}
\begin{aligned} 
s_{ij}&=\frac{1}{D}\tilde s_{ij},\\
 t_{ij}&=\frac{1}{\zeta} \de_{ij}\th_i,
\end{aligned}
\end{equation}
where $D, \zeta$ are some
non-zero numbers, satisfy the relations of $SL_2(\Z)$: $(st)^3=s^2,
s^4=\id$. These matrices are naturally interpreted as matrices of some
operators $s, t$ acting on $K(\C)$; for example, the operator
$s=\frac{1}{D}\tilde s$ where 
\begin{equation}\label{sv}
\tilde s[V]=\sum \tilde s_{VV_j}[V_j]
\end{equation}
where $[V]$ is the class in $K$ of $V\in \C$. 

We also note that the numbers $D, \zeta$ appearing in \eqref{st} are
determined uniquely up to a simultaneous change of sign. The number
$D=(s_{00})^{-1}$ is sometimes called the {\em rank} of $\C$. If $\C$
is Hermitian category, it is possible to choose $D$ to be a positive
real number. In modular tensor categories coming from conformal field
theory, the number $\zeta$ is given by $\zeta=e^{2\pi\i c/24}$, where
$c$ is the (Virasoro) central charge of the theory.

In this section we assume that $\C$ is a modular tensor category and
$A$ is a rigid $\C$-algebra, which satisfies $\th_A=\id$; by
\thref{trigidsemisimple}, this implies that $A$ is semisimple.  We
denote isomorphism classes of simple objects in $\Rep A$ by $X_\pi,
\pi\in\Pi$, and let $K(A)$ be the fusion algebra of $\Rep A$.
Similarly, set of simple objects in $\Rep^0 A$ is $\Pi^0\subset \Pi$
(see \thref{tssrep0}) and the fusion algebra of $\Rep^0 A$ is
$K^0(A)\subset K(A)$. We will denote by $P\colon K(A)\to K^0(A)$ the
projection operator: $P([X_\pi])=[X_\pi]$ if $\pi\in \Pi^0$ and
$P([X_\pi])=0$ otherwise. 
 
Define operator $\tilde s^A\colon K^0(A)\to K^0(A)$ in
the same way as for $\C$ but replacing $V_j$ by
$X_\pi$ and using rigidity morphisms in $\Rep^0A$ rather than in $\C$.

\begin{theorem}\label{tmodular} Let  $G\colon K^0(A)\to K$ be the map
  induced by the functor $G$ from \thref{tfunctors}, and let
  $F^0\colon K\to K^0(A)$ be the composition $PF$, where $P\colon
  K(A)\to K^0(A)$ is the projection operator defined above. Then $G,
  F^0$  commute with the action of $\tilde s, \tilde t$ up
to a constant:
\begin{alignat*}{2}
G(\dim A) \tilde s^A&=\tilde s G, \qquad & 
                  F^0 \tilde s&= (\dim A) \tilde s^AF^0  \\ 
G\tilde t^A&=\tilde t G,          & 
                  F^0\tilde t &=\tilde t^A F^0.
\end{alignat*} 
\end{theorem}

To prove this theorem, we will need several technical lemmas. 

\begin{lemma}\label{lstilde}
For $X,Y\in \Rep^0 A$, the number $\tilde s^A_{XY}$ is given by 

\begin{equation}
\tilde s^A_{XY}=\frac{1}{(\dim A)^2}\fig{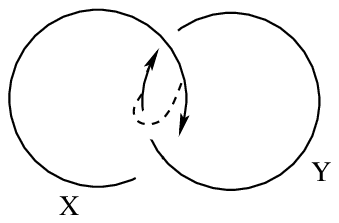}
\end{equation}

\end{lemma}
\begin{proof}
Recalling the definition of rigidity isomorphisms in $\Rep A$ and
isomorphisms $A\tta A\isoto A$, we see that  $\tilde s^A_{XY}$ is
given by 
\begin{equation}
\frac{1}{(\dim A)^2}\fig{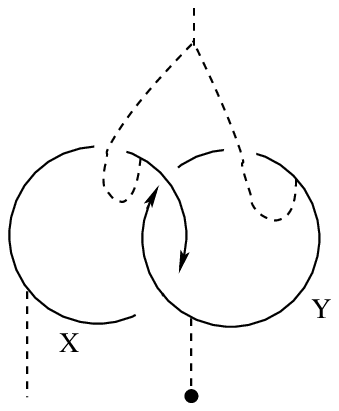}=\tilde s^A_{XY}\id_A
\end{equation}

Restricting both sides to $\one\subset A$, we get 
\begin{equation}
\tilde s^A_{XY}=\frac{1}{(\dim A)^2}\fig{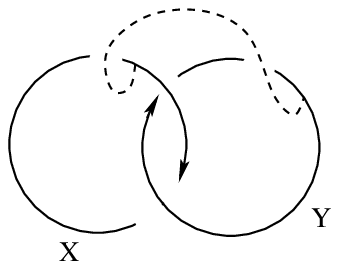}
\end{equation}
which is easily seen to be equivalent to the statement of the lemma.
\end{proof}

\begin{lemma}\label{lproj}
Let $X_\pi\in \Rep A$ be simple. Define $P_\pi\colon X_\pi\to X_\pi$
by 
\begin{equation}
P_\pi= \frac{1}{\dim A}\fig{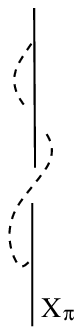}=\frac{1}{\dim A}\fig{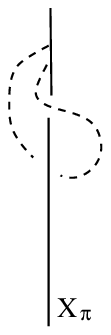}
\end{equation}
Then 
\begin{equation}
P_\pi=\begin{cases} \id_{X_\pi}, & \text{ if }X_\pi\in \Rep^0 A\\
                    0,          &\text{ otherwise}
      \end{cases}
\end{equation}      
\end{lemma}    
\begin{proof}
  If $X_\pi\in \Rep^0$, then the statement immediately follows from
  \leref{lnoose}. Thus, let us assume that $X_\pi\notin \Rep^0 A$ and
  prove that in this case, $P_\pi=0$.

First, note that the composition  $\th^{-1}_\pi P_\pi$ can be rewritten
as shown in \firef{fthP}. 

\begin{figure}[h]
\begin{equation*}
\fig{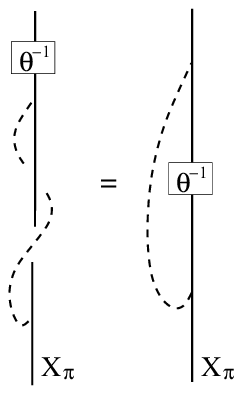}
\end{equation*}
\caption{Presentation of $\th^{-1}_\pi P_\pi$}\label{fthP}
\end{figure}

From this presentation one easily sees that $\th^{-1}_\pi P_\pi$ is a
morphism of $A$-modules; since $X_\pi$ is simple, this implies 
\begin{equation}\label{thp=c}
\th^{-1}_\pi P_\pi=c_\pi \id
\end{equation}
for some $c_\pi \in \Cset$. 

Next, let us calculate $P_\pi^2$:
\begin{equation}
\begin{aligned}
P_\pi^2&=\frac{1}{(\dim A)^2}\fig{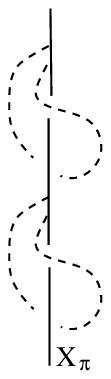}
       =\frac{1}{(\dim A)^2}\fig{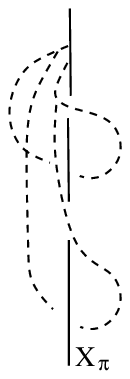}\\
       &=\frac{1}{(\dim A)^2}\fig{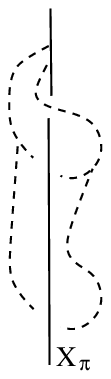}
       =\frac{1}{\dim A}\fig{proj2.eps}\\
       &=P_\pi.
\end{aligned}
\end{equation}
Thus, $P_\pi$ is a projector. On the other hand, it follows from
\eqref{thp=c} that $ P_\pi=c_\pi \th_{\pi}$. Combining these two
results, we get
$c_\pi^2\th_\pi=c_\pi$. If we assume that $c_\pi\ne 0$, then this
implies that $\th_\pi=c^{-1}_\pi$; by \thref{tssrep0}, this is impossible if
$X_\pi\notin \Rep^0 A$. Thus, $c_\pi=0$. 
\end{proof}

\begin{lemma}\label{lsxy}
For $X\in \Rep^0 A, Y\in \Rep A$, one has 
\begin{equation}
\langle\tilde s^A(X), Y\rangle
 =\langle\tilde s^A(P(Y)), X\rangle
 =\frac{1}{(\dim A)^2} \fig{stilde.eps}
\end{equation}
where $P\colon K(A)\to K^0(A)$ is as in \thref{tmodular}.
\end{lemma}
\begin{proof}
Since both sides are linear in $Y$ it suffices to prove this
formula when  $Y$ is simple. If $Y\in \Rep^0 A$, the
statement immediately follows from \leref{lstilde}. Thus, we only need
to prove that if $Y$ is simple, $Y\notin \Rep^0A$, then the
right-hand side is zero. To prove this, let $C\in\Cset$ be defined by 
$$
C=\fig{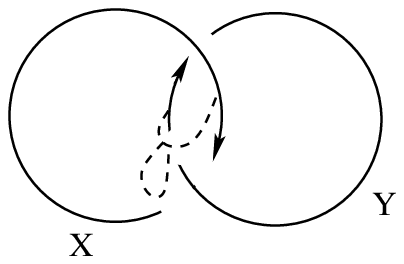}
$$

On one hand, it easily follows from \leref{lnoose} that 
$$
C=\dim A \fig{stilde.eps}
$$
On the other hand, we can deform the figure defining $C$ as shown
below
\begin{equation}
\begin{aligned}
C&=\fig{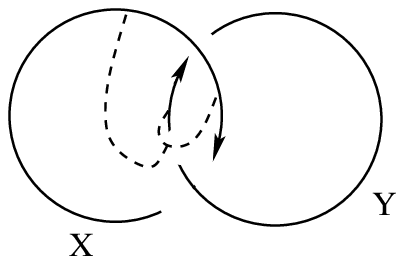}=\fig{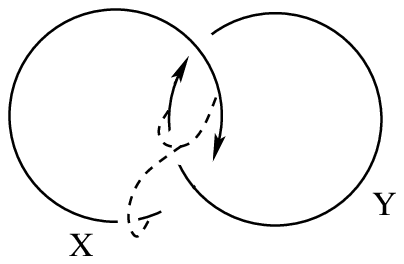}\\
              &=\fig{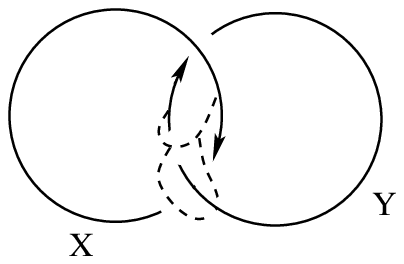}=\fig{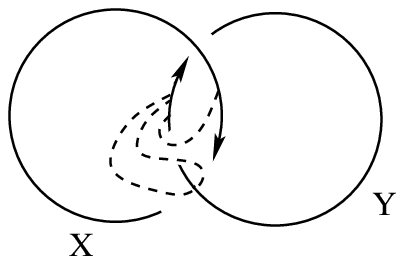}
\end{aligned}
\end{equation}

By \leref{lproj}, this implies $C=0$.
\end{proof}

Now we are ready to prove \thref{tmodular}.
\begin{proof}[Proof of \thref{tmodular}]
  Proof for $\tilde t$ is obvious from the definition. As for $\tilde
  s$, it suffices to prove that $(\dim A)\langle G\tilde s^A (X), V\rangle
  =\langle \tilde s (G (X)), V\rangle$ for any $X\in \Rep^0 A, V\in
  \C$. Using adjointness of $G$ and $F$, this reduces to
$$
(\dim A)\langle \tilde s^A (X), F(V)\rangle=\tilde s_{G(X),V}.
$$
 (note that $F(V)\in \Rep A$, but in general, not in $\Rep^0 A$).
Using \leref{lsxy} and definition of $F(V)$,  this can be rewritten as
the following identity of figures:

\begin{equation}\label{keyeq} 
\frac{1}{\dim A}\fig{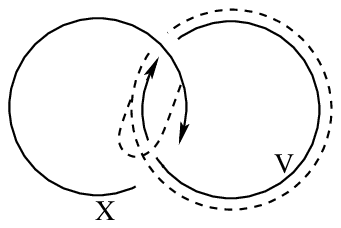}=\fig{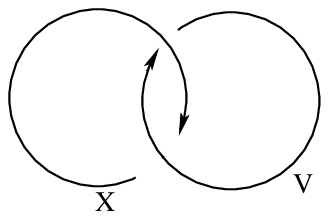}
\end{equation}
which can be proved by rewriting the graph in left hand side as shown in
\firef{fsxfv} and using \leref{lproj}. 

\begin{figure}[ht] 
$$
\fig{sxfv}
 =\fig{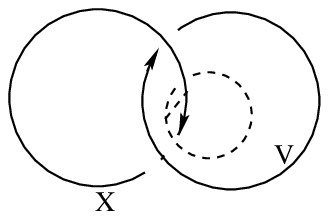}
 =\fig{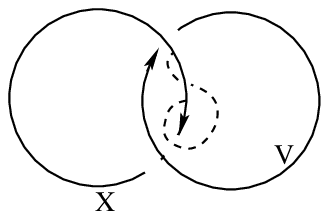}
$$
\caption{}\label{fsxfv}
\end{figure}

Similarly, the identity involving  $F^0$ is equivalent to 
$$
\langle(\dim A) \tilde s^A(F^0(V)), X\rangle=
\langle\tilde s(V), G(X)\rangle
$$
which is also equivalent to \eqref{keyeq}.  This completes the
proof of \thref{tmodular}.
\end{proof} 

This theorem implies the following important result. 

\begin{theorem}\label{tmodular2}
  If $\C$ is a modular category, $A$ is a rigid $\C$-algebra,
  $\th_A=\id$, then $\Rep^0 A$ is modular and the numbers $D, \zeta$
  appearing in \eqref{st} for $\Rep^0 A$ are related with the
  corresponding numbers for $\C$ by
\begin{equation}\label{dzeta}
\begin{aligned}
D(\Rep^0A)&=\frac{D(\C)}{\dim A}\\
\zeta(\Rep^0 A)&=\zeta(\C).
\end{aligned} 
\end{equation}
Also, the maps $G\colon K(\Rep^0 A)\to K(\C), F^0\colon K(\C)\to
K(\Rep^0 A)$ commute with operators $s, t$.
\end{theorem}
\begin{proof}
The proof is based on the following lemma. 
\begin{lemma}
Let $\A$ be a semisimple rigid braided balanced category over $\Cset$, with
finitely many isomorphism classes of simple objects. Then $\A$ is
modular iff the matrix $\tilde s$, defined by \firef{fsij}, satisfies 
\begin{equation}\label{s^2}
\tilde s^2\one =c\one
\end{equation}
for some $c\in \Cset, c\ne 0$. 
\end{lemma}
This lemma is not new; however, for the sake of completeness, we
include its proof below. 

Thus, to prove that $\Rep^0A$ is modular, it suffices to prove
$(\tilde s^A)^2 A=cA$ for some $c\ne 0$. But by \thref{tmodular},
$\tilde s^A$ commutes with $F^0$ up to a constant; thus,
\begin{align*}
(\tilde s^A)^2 A&=(\tilde s^A)^2 F^0(\one)
       =\frac{1}{(\dim A)^2}F^0(\tilde s^2\one)
       =\frac{1}{(\dim A)^2}F^0(c\one)\\
    &=\frac{c}{(\dim A)^2}A. 
\end{align*}
Thus, $\Rep^0A$ is modular; all other statements of the theorem
immediately follow from \thref{tmodular}.

\begin{proof}[Proof of the lemma]
If $\A$ is modular, the statement is well known and in fact $c=D^2$
(see, e.g., \cite{BK}). Thus, let us assume that \eqref{s^2} holds and
deduce from it non-degeneracy of $\tilde s$. 

First, note that $\tilde s\one =\dd=\sum d_i V_i$, where $V_i$ are simple
objects in $\A$ and $d_i=\dim_\A V_i$. Thus, \eqref{s^2} implies
$\langle \tilde s\dd, V_i\rangle =c\de_{i,0}$ which can be rewritten
as 
\begin{equation}\label{s^2a}
\sum_{j}d_j \fig{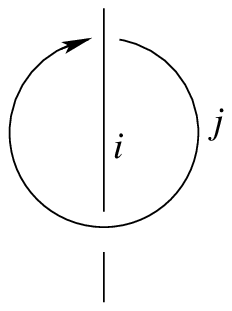}=c\de_{i,0}\fig{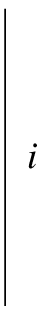}
\end{equation}

Let us now choose some $i,k\in I$ and let 
$$
C_{ijk}=\fig{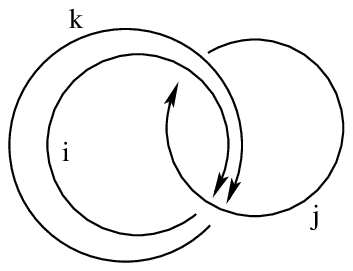}
$$
On one hand, it is easy to show using the definition of $\tilde s$ that 
$$
C_{ijk}=\frac{\tilde s_{ij}\tilde s_{jk}}{d_j}
$$
(see, e.g., \cite[Lemma 3.1.4]{BK}) and thus, 
$$
\sum_j d_jC_{ijk}=\sum_j\tilde s_{ij}\tilde s_{jk}=(\tilde s ^2)_{ik}.
$$
On the other hand, decomposing $V^*_i\ttt V^*_k$ in a direct sum of
irreducibles and using  \eqref{s^2a}, we get 
$$
\sum_j d_jC_{ijk}=c\langle V^*_i\ttt V^*_k,\one\rangle =c\de_{ik^*}
$$
which is a non-singular matrix. Therefore, $\tilde s^2$ is
non-singular and thus $\A$ is modular. This completes the proof of
the lemma and thus of \thref{tmodular2}.
\end{proof}
\renewcommand{\qed}{}
\end{proof}
\begin{remark}
For modular tensor categories coming from conformal field theory, the
identity $\zeta(\C)=\zeta(\Rep^0 A)$ can be interpreted as stating that
an extended CFT has the same central charge as the original CFT,
which, of course, should be expected. 
\end{remark}

\section{Vertex operator algebras}\label{svoa}
In this section, we give the example which was one of our main
motivations for this work. Detailed proofs of the results given here
will appear in a separate paper \cite{voaproof}; here we only outline
the main ideas. We should also note that relation between extensions
of vertex operator algebras and algebras in a category discussed here
was independently found by A.~ Wassermann \cite{was}, who discussed
his work with the second author (Ostrik) during his visit to MSRI in
November 2000. At this time, we were finishing the first draft of this
paper.

We assume that the reader is familiar with the notion of a vertex
operator algebra (VOA); a review and list of references can be found,
e.g., in \cite{Frenkel}, \cite{FHL}. To avoid ambiguity, we mention
that we include the Virasoro element $\omega$ and $\Z$-grading in the
definition of a VOA.  Similarly, when talking about modules over a
VOA, we always assume that $L_0$ acts semisimply with
finite-dimensional eigenspaces (this automatically gives
$\Cset$-grading on a module). We only consider finite length modules.

 Let $\V$ be a vertex operator algebra which
is nice enough so that the following properties are satisfied:

\begin{enumerate}
  
\item For every simple $\V$-module $M$, its conformal dimension (i.e.
  lowest eigenvalue of $L_0$) is real, $\ge 0$, with equality only for
  $M=\V$, in which case $\dim \V_0=1$. 

\item The category of representations of $\V$ is semisimple, with only
  finitely many simple objects, and all spaces of conformal blocks
  (i.e., intertwining operators between tensor products of
  representations) are finite-dimensional. Also, $\V$ is simple as a
  $\V$-module. 

\item The category $\C$ of $\V$-modules is a rigid braided tensor
  category. 
\end{enumerate}
  
The first condition is technical; we will only need it to ensure
uniqueness of vacuum vector (see proof of \thref{thvoa} below). The
most important condition is the last one, which deserves detailed
discussion. 

There are at least two ways to define the tensor product (usually
called the fusion tensor product) structure on the category of
$\V$-modules, both originating in the pioneering work of Moore and
Seiberg.  The first construction, developed in a series of papers of
Huang and Lepowsky \cite{tensor12}, \cite{tensor3}, \cite{tensor4}, is
based on defining the tensor product via the space of intertwining
operators. The second approach uses the vector spaces of coinvariants
(see \cite{Frenkel}) which should give a modular functor, and then
using this modular functor to define the structure of a tensor
category (see \cite{BK}). This shows that for every VOA appearing in
conformal field theory the category of modules has a structure of a
rigid braided tensor category.  In fact, such a VOA has to satisfy a
stronger restriction:
\begin{enumerate}
\item[($3'$)]
The category of $\V$-modules is a modular tensor category. 
\end{enumerate}
Indeed, it follows from the axioms of a rational conformal field
theory that the spaces of conformal blocks for such a VOA form a
modular functor, and it is known that a modular functor allows one to
define a structure of a modular tensor category (see, e.g.,
\cite{BK}).

These two approaches should give equivalent results; unfortunately, to
the best of our knowledge, details of this equivalence are not
available in the literature. In what follows, we will use the first
approach, i.e. use the definition of the tensor structure given by
Huang and Lepowsky.

In both approaches, it is relatively easy to give the definition of
the tensor product, but it is extremely difficult to check that for a
given VOA this tensor product is well-defined and defines a structure
of a rigid balanced braided category (see \cite{tensor4} for the list
of conditions that need to be checked).  So far, this has only been
checked in very few examples.\footnote{Of course, as mentioned above,
  this should hold for any VOA that comes from a rational conformal
  field theory, but this does not help much: axioms of RCFT are even
  more difficult to check.}

Most important example of a VOA for which conditions (1)--(3) have
been checked is the VOA coming from an  affine Lie algebra at positive
integer level, discussed below.

\begin{example}\label{exwzw}
  Let $\g$ be a simple Lie algebra, $\ghat$ --- corresponding affine
  Lie algebra, and $k$ --- a non-negative integer (level). Let
  $L_{0,k}$ be the integrable $\ghat$ module of level $k$ with highest
  weight $0$ (the vacuum module). Then it is known that it has a
  canonical structure of a VOA; we will denote this VOA by $\V(\g,
  k)$. This VOA satisfies requirements (1)--($3'$) and thus, its
  category of representations $\C(\g,k)$ is modular (see \cite{HL},
  \cite{BK}). As an abelian category, $\C(\g,k)$ is just the category
  of integrable $\ghat$ modules of level $k$. It is also known (see
  \cite{Finkelberg}) that $\C(\g,k)$ is equivalent (as modular
  category) to the ``semisimple part'' of the category of
  representations of the quantum group $U_q\g$ with
  $q=e^{\pi\i/m(k+h^\vee)}$, where $h^\vee$ is the dual Coxeter number
  and $m=1$ for simply-laced algebras, $m=2$ for $B_n, C_n, F_4$ and
  $m=3$ for $G_2$.
\end{example}

Let $\V\subset \V_e$ be a subalgebra (in the sense of VOA's). Assume
in addition that $\V_e$ is finite length as a module over $\V$. Then
we will call $\V_e$ an {\em extension} of $\V$.

\begin{theorem}\label{thvoa}
  Let $\V$ be a VOA satisfying {\rm{(1)--(3)}} above, and let $\C$ be
  the category of $\V$--modules. Then the following two notions are
  equivalent:
\begin{enumerate}
\item
 An extension $\V\subset \V_e$, where $\V_e$ is also a VOA satisfying
 properties {\rm{(1)--(3)}} above 

\item A rigid $\C$-algebra $A$ with $\th_A=1$
\end{enumerate}

  Under this correspondence, category of
   $\V_e$-modules is identified with $\Rep^0A$.
\end{theorem}
\begin{proof}
  We give a sketch of the proof; details will appear in the
  forthcoming paper \cite{voaproof}.
  
  If $\V_e$ is a VOA, then for every $v\in \V_e$ we have the vertex
  operator $Y(v,z)\colon \V_e\to \V_e[[z,z^{-1}]]$. Restricting it to
  $v\in \V$, we get a structure of a $\V$-module on $\V_e$. It is
  immediate from the definitions that the map $Y(\cdot, z)$ is an
  intertwining operator of the type $\binom{\V_e}{\V_e\ \V_e}$ and
  thus gives a morphism of $\V$-modules $Y\colon \V_e\stackrel{\bf
    .}{\otimes}\V_e\to \V_e$, where $\stackrel{\bf .}{\otimes}$ is the
  ``fusion'' tensor product. It follows from the usual commutativity
  and associativity axioms for a VOA (see \cite{FHL}) that $Y$ defines
  a structure of a commutative and associative algebra on $\V_e$.
  Existence and uniqueness of unit follow from existence and
  uniqueness of the vacuum vector in a VOA (see condition (1) above).
  Condition $\th_A=1$ follows from the fact that eigenvalues of $L_0$
  on $\V_e$ are integer. A straightforward check shows that the
  arguments above can be reversed and that the category of
  representations of $\V_e$ as a VOA coincides with $\Rep^0 A$.
\end{proof}

One of the general ways to construct extensions of the VOA $\V(\g,k)$
is by using the notion of conformal embedding (note, however, that not
all extensions can be obtained in this way). Let $\g\subset \g'$ be an
embedding of Lie algebras; then it defines an embedding of affine Lie
algebras $\ghat\subset \ghat'$. This embedding doesn't preserve the
level --- a pullback of a $\ghat'$ module of level $k'$ will be a
module of level $k=x_e k'$ for some integer $x_e$; we will
symbolically write $(\ghat)_k\subset (\ghat')_{k'}$. It defines an
embedding $\V(\g, k)\subset \V(\g', k')$ which preserves the operator
product expansion (i.e., the algebra structure in $\V$) but in general
not the Virasoro element. In some special cases, however, such an
embedding preserves the Virasoro element as well and therefore defines
an embedding of VOA's; they are called {\em conformal embeddings}.
In this case it is easy to show (see, e.g., \cite[Chapter 17]{CFT})
that $\V(\g', k')$ is automatically finite as $\ghat$-module, so
$\V(\g', k')$ is an extension of $\V(\g,k)$.

\begin{example}\label{confembed}
  Let $\C(\slt,k)$ be the category of integrable modules over
  $\slthat$ of level $k$. Then it is known that for $k=10$, there is a
  conformal embedding $(\slthat)_{10}\subset \widehat{sp(4)}_1$. The
  easiest  way to describe this embedding is to note that the
  irreducible 4-dimensional representation of $\slt$ has an invariant
  non-degenerate skew-symmetric form, which gives an embedding
  $\slt\subset sp(4)$.
  
  The decomposition of $\V(sp(4), 1)$ as $\V(\slt, 10)$ module is
  given by $\V=L_{0,10}\oplus L_{6,10}$ (see \cite[Chapter 17]{CFT}).
  Thus, this shows that the object $A=L_{0,10}\oplus L_{6,10}\in
  \C(\slt, 10)$ has a structure of a rigid $\C$-algebra (later
  we will show that such a structure is unique).
  
  Similarly, for $k=28$ there exists a conformal embedding
  $(\slthat)_{28}\subset (\widehat{G}_2)_1$; the decomposition of
  $\V(G_2,1)$ as $\V(\slt,28)$ module is given by $\V=L_{0,28}\oplus
  L_{10, 28}\oplus L_{18, 28}\oplus L_{28,28}$.
\end{example}

\begin{remark}
The use of conformal embeddings to produce extensions of chiral
algebras is, of course, well known in physics literature. Conformal
embeddings can also be used in the subfactor theory --- see \cite{X}
and references therein. 
\end{remark}
\begin{remark} 
  It is very easy to prove rigorously that the conformal embedding
  $\ghat_k\subset \ghat'_{k'}$ determines a structure of a rigid
  $\C-$algebra on the vacuum module $V=L_{0,k'}$ over $\ghat'$
  considered as a module over $\ghat$, even without referring to the
  more general \thref{thvoa}. Let $\otimes_\g, \otimes_{\g'}$ denote
  the fusion product over $\ghat, \ghat'$ respectively. this fusion
  tensor product can be defined in terms of coinvariants for the
  action of the algebra of rational $\g$ (respectively, $\g'$) valued
  functions.  Namely, consider the rational curve ${\mathbb P}^1$ with
  3 marked points and the representations $V, V, V^*$ assigned to
  these points. The spaces of homomorphisms $V\otimes_{\g} V\to V,
  V\otimes{\g'}$ are isomorphic to the dual of the spaces of
  coinvariants for $\g, \g'$, see \cite{KL, Finkelberg}. Also, the space of
  coinvariants for $\g'$ is canonically isomorphic to $\C$. By
  definition, we have a surjection $Coinv(\g)\surjto Coinv(\g')$.
  Thus, we have an embedding $\Hom (V\otimes_{\g'} V, V)\subset
  \Hom(V\otimes_{\g} V, V)$. Thus, the canonical morphism
  $V\otimes_{\g'} V\to  V$ defines a morphism $m:
  V\otimes_h V\to V$. Let us prove that this morphism defines an
  associative multiplication, that is $m(m\otimes_{\g} id)=m(id\otimes_{\g}
  m)$. Both sides of this equality are represented by some
  coinvariants (for ${\mathbb P}^1$ with 4 marked points and the
  representations $V, V, V, V^*$ assigned to these points) and by the
  construction of the fusion product these coinvariants actually come
  from the coinvariants over $g$.  But the space of the coinvariants
  over $\g'$ is one dimensional since $V\otimes_{g'} V\otimes_{\g'}V\simeq V$
  and hence the LHS and the RHS are proportional. To compute the
  proportionality coefficient it is enough to note that $m$ restricts
  nontrivially on $V_0\subset V$ where $V_0$ is the vacuum module over
  $\ghat$. The proof of commutativity is completely analogous.
  Finally one can use the coinvariant above to identify $V$ and $V^*$
  as $\ghat-$modules what implies that $m\colon V\otimes_h V\to V\to V_0$
  can be used to identify $V$ and $V^*$ as $\ghat-$modules that is
  $V$ is rigid $\C-$algebra over $\ghat$.
\end{remark}
\section{ADE classification for $\U$}
In this section, we apply the general formalism developed above in a
special case: when $\C$ is the semisimple part of the category of
representations of $\U$ with $q=e^{\pi\i/l}, l\leq 2$ as defined by
Andersen et al \cite{AP}. We assume that the reader is familiar with the
definition and main properties of categories of representations of
quantum groups at roots of unity; if not, we refer to the monograph
\cite{BK} for a review.

It is known that the category $\C$ is semisimple, with simple objects
$V_0, \dots, V_k$, where $k=l-2$ and $V_i$ is the usual
$(i+1)$-dimensional irreducible representation of $\U$. Its
Grothendieck ring $K$ is generated by one element, $V_1$. The quantum
dimensions are given by
$$
\dim_\C V_n=[n+1]:=\frac{q^{n+1}-q^{-(n+1)}}{q-q^{-1}}
$$
which in particular implies that for any  non-zero object $V$,
\begin{equation}\label{d>1}
\dim_\C V\ge 1.
\end{equation}
It is also known that this category is modular and that the universal
twist $\th$ is given by 
\begin{equation}\label{thi}
\th_n:=\th_{V_n}=q^{n(n+2)/2}=e^{2\pi\i\frac{n(n+2)}{4(k+2)}}.
\end{equation}

Finally, we note that this category is equivalent to the category
of integrable representations of affine Lie algebra $\widehat{\slt}$
of level $k=l-2$, or, equivalently, the category of representations of
the corresponding vertex operator algebra $\V(\slt, k)$ (see
\cite{Finkelberg}).

Our main goal is to classify all $\C$-algebras.


\begin{theorem} \label{tade}
  There is a correspondence between rigid  $\C$-algebras with
  $\th_A=\id$  and
  Dynkin diagrams of types $A_n, D_{2n}, E_6, E_8$ with Coxeter number
  equal to $l$. Under this correspondence, simple objects of $\Rep A$
  correspond to vertices of the Dynkin diagram, and the matrix of
  multiplication by $F(V_1)$ in $K(\Rep A)$ is $2-A$, where $A$ is the
  Cartan matrix of the Dynkin diagram. 
\end{theorem}

\begin{proof}
  Let $A$ be a rigid $\C$-algebra with $\th_A=\id$. In this case,
  $\Rep A$ is a monoidal category and, by \coref{cmodule}, a module
  category over $\C$. This implies that the Grothendieck ring $K(A)=K(\Rep
  A)$ is a module over $K(\C)$. By \thref{trigidsemisimple}, $\Rep A$
  is semisimple, so $K(A)$ has a distinguished basis (classes
  $[X_\pi]$ of simple objects) so that in this basis, multiplication
  by any $F(V), V\in \C$, has coefficients from $\Z_+$. In addition,
  this module has the following properties:

\begin{enumerate}
\item[(i)] The module $K(A)$ is indecomposable: it is impossible to
  split the set of simple objects $\Pi$ as $\Pi=\Pi'\sqcup \Pi''$ so
  that $K'=\oplus_{\Pi'}\Cset [X_\pi], K''=\oplus_{\Pi''}\Cset
  [X_\pi]$ are $K(\C)$-submodules in $K(A)$. 

  Indeed, every simple module $X_\pi$ appears with non-zero
  multiplicity in $F(V_i)\tta A=F(V_i)$ for some $V_i$. This follows
  from $\langle F(V_i),X_\pi\rangle=\langle V_i,G(X_\pi)\rangle$.

\item[(ii)] There exists a map $d\colon K(A)\to \Cset$ such that
  $d(X_\pi)\in\Rset_{>0}$ and $d(F(V)\tta X)=(\dim_\C V) d(X)$. 
  
  Indeed, it suffices to let $d(X)=\dim_{\Rep A}(X)$ and use
  \thref{tdim0}. 

\item[(iii)] There exists a symmetric bilinear form $\langle, \rangle$
  on $K(A)$ such that $\langle F(V)\tta X, Y\rangle=\langle X,
  F(V)\tta Y \rangle$ for any $V\in \C, X,Y\in\Rep A$. 
  
  Indeed, we can let $\langle X, Y\rangle=\dim\Hom_A(X,Y)$ and use
  rigidity and $F(V_i)^*\simeq F(V_i^*)\simeq F(V_i)$ (not
  canonically).
\end{enumerate}

All modules $M$ over $K(\C)$ which have properties (i)--(iii) above
were classified in \cite{EK}, where it is shown that they correspond
to finite Dynkin diagrams with loops with Coxeter number equal to $l$.
Under this correspondence, vertices of the Dynkin diagram correspond
to the elements of distinguished basis of $M$, and the matrix of
multiplication by $V_1\in \C$ is $2-A$, where $A$ is the Cartan matrix
of the Dynkin diagram.

(Dynkin diagrams with loops, in addition to the usual Dynkin diagrams,
include ``tadpole'' diagrams  $T_n$ shown in \firef{diagrL}; in
\cite{EK}, this diagram is denoted by $L_n$. By
definition, the Cartan matrix for such a diagram is the same as for
$A_n$ but with $a_{11}=1$, and the Coxeter number for $T_n$ is equal
to $2n+1$).

\begin{figure}[h]
\begin{equation*}
\fig{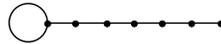}
\end{equation*}
\caption{Dynkin diagram of type $T$.}\label{diagrL}
\end{figure}

\begin{remark}  
  In an interesting note \cite{mckay2}, it was shown that the dimension
  vector $d$ can also be obtained from so-called ``semi-affine'' Dynkin
  diagrams, which give $d$ both for finite and affine
  Dynkin diagrams.
\end{remark}

Now we have to check which of these modules can actually appear as
Grothendieck ring $K(A)$ for some rigid $\C$-algebra $A$. 

First, note that if $K(A)$ is indeed the Grothendieck ring of a
rigid $\C$-algebra $A$, then not only we have a
distinguished basis $[X_\pi]$ and an inner product $\langle,\rangle$
but in fact, the distinguished basis is orthonormal with respect to
$\langle,\rangle$. This implies that the matrix of tensor product with
$F(V_1)$ is symmetric in this basis. Thus, only simply-laced Dynkin
diagrams can possibly come from $K(A)$. This leaves us with the ADET
type diagrams.

Next, we need to determine which vertex of the Dynkin diagram
corresponds to the unit object, i.e. to $A$ itself. 

\begin{lemma}\label{lleg}
If $A$ is a rigid $\C$-algebra, then $A$ corresponds to the
end of the longest leg of the corresponding Dynkin diagram. 
\end{lemma}
\begin{remark}
By an ``end'' we mean a vertex which is connected to exactly one
vertex; in particular, the vertex with a loop in the diagram of type
$T$ is not considered an end vertex.
\end{remark}

\begin{proof}
  Let $X\in \Rep A$ be the object corresponding to one of the
  ends of legs of the Dynkin diagram. Then $F(V_1)\tta X$ is simple.
  Since in a rigid category, tensor product of non-zero objects is
  always non-zero, this implies that $F(V_1)\simeq F(V_1)\tta A$ is
  simple. Thus, $A$ is connected to exactly one vertex, which
  means that $A$ itself is an end of one of the legs. 
 
  To prove that $A$ is the end of the longest leg, note that if $X$ is
  an end of the leg of length $m$ (that is, consisting of $m$ edges),
  then $F(V_1)\tta X,\dots, F(V_m\tta X)$ are simple but
  $F(V_{m+1})\tta X$  is not. This implies that $F(V_i)=F(V_i)\tta A,i=1,
  \dots, m$ are simple, which means that the leg containing $A$
  has length at least $m$. 
\end{proof}




This determines the vertex corresponding to $A$ uniquely up to an
automorphism of the Dynkin diagram. 

Once we know the vertex corresponding to $A$, we know the class of
$F(V_1)$ in $K(A)$; since $F$ is a tensor functor and $V_1$ generates
$K$, this uniquely determines the map $F$ at the level of Grothendieck
rings, and thus, the adjoint map $G\colon K(A)\to K$. In other words,
we can write for each vertex of the Dynkin diagram the structure of
the corresponding object $X_\pi$ as an object of $\C$. In particular,
this gives decomposition of $A$ itself as an object of $\C$.

Doing this explicitly for diagrams $A_n, D_n, E_n, T_n$ gives the
answer shown in Table~\ref{table1} (no, it was not found using a
computer --- it is done easily by hand), which agrees with the one
given in Cappelli-Itzykson-Zuber classification.
\begin{table}[h]
\caption{Algebra $A$ for various Dynkin diagrams}\label{table1}
\begin{tabular}{c|c|c}
Diagram& $k=h-2$ & $A$\\
\hline  
$A_n$ & $n-1$ & $V_0$ \\
$D_n$ & $2n-4$ & $V_0+V_k$ \\
$T_n$ & $2n-1$ & $V_0+V_k$ \\
$E_6$ & $10$ & $V_0+V_6$ \\
$E_7$ & $16$ & $V_0+V_8+V_{16}$ \\
$E_8$ & $28$ & $V_0+V_{10}+V_{18}+V_{28}$ \\
\end{tabular}
\end{table}

Next step is to find which of the possible $A$ given in this table do
have a structure of a $\C$-algebra. 

{\bf Type $A$}: in this case, $A=\one$ obviously has a unique
structure of commutative associative algebra, and $\Rep A=\C$. 

{\bf Type $D$}.  Let us introduce the notation 
\begin{equation}\label{de}
\de=V_k.
\end{equation}
It easily follows from explicit formulas that $\dim_\C \de=1$ and
$\de\ttt V_n\simeq V_{k-n}$; in particular, $\de\ttt\de\simeq \one$. 

\begin{theorem}\label{tdtype}
The object $A=\one\oplus\de$ in $\C$ has a structure of a rigid 
$\C$-algebra iff $4|k$. In this case, the
structure of an algebra is unique up to isomorphism, and this algebra
satisfies $\th_A=\id$. 
\end{theorem}
\begin{proof}
  Let $\mu$ be the multiplication map $\mu\colon
  (\one\oplus\de)\ttt(\one\oplus\de)\to (\one\oplus\de)$. All
  components of such a map are uniquely determined by the unit axiom,
  except for $\mu_{\de\de}\colon \de\ttt\de\to \one$. Since
  $\de\ttt\de\simeq \one$, such a map is unique up to a constant.
  Rigidity implies that $\mu_{\de\de}\ne 0$. This proves
  uniqueness. 

  To check existence, fix some non-zero $\mu_{\de\de}$. Then
  associativity and commutativity are equivalent to 
\begin{equation}
\begin{aligned}
\mu_{\de\de}\circ(\id \ttt\mu_{\de\de})
 &=\mu_{\de\de}\circ(\id \ttt\mu_{\de\de})\colon \de\ttt\de\ttt\de\to\de\\
\mu_{\de\de}\circ R_{\de\de}&=\mu_{\de\de}.
\end{aligned}
\end{equation}
To check the second equation, we use the following lemma

\begin{lemma}\label{lfR}
For generic values of $q$, let $f\colon V_{a}\otimes V_{a}\to V_{2b}$
be a nonzero homomorphism. Then 
$$
f\circ R_{V_a V_a}=(-1)^{a-b}\th_a^{-1} (\th_{2b})^{1/2}f
$$
where $\th_{a}=q^{a(a+2)/2}$ is the universal twist and
$\th_{2b}^{1/2}=q^{2b(2b+2)/4}$. 
\end{lemma}

To prove this lemma, note that it immediately follows from balancing
axiom in $\C$ that $f\circ R^2=\th_a^{-2}\th_{2b}$, which gives the
formula above up to a sign. To find the sign, it suffices to let
$q=1$.

Since this formula works for generic values of $q$, it should also be
valid for $q$ being a root of unity. In particular, applying this
lemma to 
$q=e^{\pi\i/(k+2)}$ and $\mu_{\de\de}\colon \de\ttt\de\to\de$, we get 
$$
\mu_{\de\de}\circ R_{\de\de}=(-1)^k \th^{-1}_{\de} \mu_{\de\de}.
$$
We have $\th_\de =q^{\frac{k(k+2)}{2}}=e^{\pi \i
  k(k+2)/2(k+2)}=e^{2\pi \i k/4}=\i^k$. Thus,
$(-1)^k\th^{-1}_\de=\i^k$ is equal to one iff $k$ is divisible by 4.
Therefore, the map $\mu$ is commutative iff $k=4m$.

To check associativity, note that both sides are equal up to a
constant (since $\dim \Hom(\de^{\ttt 3}, \de)=1$); to find the
constant, take composition of  both sides with $(i_\de)\ttt \id\colon
\de\to \de^{\ttt 3}$ and use $\dim_\C \de=1$.
\end{proof}

The category of representations of this algebra is described in detail
in \seref{sdtype}. It follows from the analysis there that the structure
of $K(A)$ as $K(\C)$-module is described by the diagram $D_{2m+2}$.

{\bf Type $T$}.

The diagram $T_n$ can not appear as $K(A)$ for a
commutative associative algebra $A$. Indeed, in this case $A$ must be
isomorphic to $V_0\oplus V_k$, but it was proved in \thref{tdtype} that
there is at most  one structure of a rigid $\C$-algebra on this
object, and if it exists, $K(A)$ is described by $D_{n}$, not $T_n$. 

{\bf Type $E_7$}.

This diagram can not appear as $K(A)$ for a commutative associative
algebra $A$. Indeed, in this case the table gives $A=V_0\oplus
V_8\oplus V_{16}=(\one\oplus\de)\oplus V_8$. Obviously,
$A'=\one\oplus\de$ is a subalgebra in $A$, and multiplication on $A$
defines a structure of $A'$-module on $V_8$ and morphism of
$A'$-modules $V_8\ttt V_8\to (\one\oplus\de)$. By rigidity, this morphism
is non-zero, which also implies that the restriction of $\mu$ to
$V_8\ttt V_8\to \de$ is non-zero. But it immediately follows from
\leref{lfR} that such a morphism can not be symmetric.

{\bf Type $E_6$}. 

In this case, there is a unique up to isomorphism $\C$-algebra
structure on $V_0\oplus V_6$. Existence follows from the discussion of
the previous section and existence of a conformal embedding of affine
Lie algebras $(\slthat)_{10}\injto \widehat {sp(4)}_{1}$ (see
\exref{confembed}).  To prove uniqueness, note that the only
non-trivial components of the multiplication map $\mu$ are $\mu'\colon
V_6\ttt V_6\to \one$, $\mu''\colon V_6\ttt V_6\to V_6$. Both of them
are unique up to a constant factor. We can fix some non-zero morphisms
\begin{align*}
e\colon &V_6\ttt V_6\to\one,\\ 
f\colon &V_6\ttt V_6\to V_6.
\end{align*}
Then $\mu'=\al e, \mu''=\be f$ for some $\al, \be\in \Cset$. It
follows from rigidity that $\al\ne 0$. Using isomorphism of
$\C$-algebras $\ph\colon (\one\oplus V_6)\to (\one\oplus V_6)$ given
by $\ph|_\one=\id, \ph|_{V_6}=\al^{1/2}\id$, we see that without loss
of generality we can assume $\al=1$, so $\mu|_{V_6\ttt V_6}=e+\be f$.
Condition that $\mu$ be associative gives the following quadratic
equation on $\be$:
$$
\be^2\Ph_1=\Ph_2
$$
where $\Ph_1, \Ph_2$ are morphisms $V_{6}^{\ttt 3}\to V_6$ given by   
\begin{align*}
\Ph_1&= f\circ (\id \ttt f)-f\circ (f\ttt \id)\\
\Ph_2&=e\ttt \id-\id \ttt e.
\end{align*}

It is easy to see that $\Phi_2\ne 0$, so the equation
$\be^2\Ph_1=\Ph_2$ is non-trivial. Thus, such an equation may either
have no solutions at all or have exactly two solutions differing by
sign: $\be=\pm \be_0$. These two solutions actually would give
isomorphic algebras: the map $\ph\colon\one\oplus\de\to\one\oplus\de$ given
by $\ph|_{\one}=1, \ph|_{\de}=-1$ gives the isomorphism.

{\bf Type $E_8$.}

In this case, there again exists a unique structure of a rigid
$\C$-algebra on $A=V_0\oplus V_{10}\oplus V_{18}\oplus V_{28}$.
Existence follows from existence of conformal embedding
$(\slthat)_{28}\subset (\widehat{G}_2)_1$ (see \exref{confembed}). To
prove uniqueness, let $A'\subset A$ be the subalgebra generated (as a
$\C$-algebra) by $V_0\oplus V_{10}$. Let $\V_e$ be the vertex operator
algebra corresponding to $A'$; by results of \seref{svoa}, it is an
extension of the VOA $\V=\V(\slt, 28)$. From the definition, $\V_e$ is
generated as a VOA by $\V$ and $L_{0,10}$. Since $L_{0,10}$ is an
irreducible $\slthat$ module, it is generated (as $\slthat$ module) by
its lowest degree component (degree stands for homogeneous degree,
i.e. eigenvalue of $L_0$). This lowest degree is equal to
$\Delta_{10}=\frac{10(10+2)/2}{2(28+2)}=1$.

Since it is well known that $\V(\g, k)$ is generated as a VOA by its
degree one component $\V[1]\simeq \g$, we see that $\V_e$ is generated
as a VOA by $\V[1]\oplus L_{10,28}[1]$. It is also easy to check that
conformal dimensions (i.e., lowest eigenvalues of $L_0$) for
$L_{18,28}$ and $L_{28,28}$ are greater than one, so
$\V_e[1]=\V[1]\oplus L_{10, 28}[1]\simeq \slt\oplus L_{10}$, where
$L_{10}=L_{10,28}[1]$ is an irreducible $\slt$-module with highest
weight 10.

By \cite[Section 2.6]{Kac}, if $\V_e[0]=\Cset, \V_e[n]=0$ for $n<0$
and $\V_e$ is generated as VOA by $\V_e[1]$, then $\V_e[1]=\g$ is a
Lie algebra with an invariant bilinear form, and $\V_e$ is naturally a
module over $\ghat$; moreover, $\V_e$ is a quotient of the Weyl module
$V^{\g}_{0,k}$ over $\ghat$ for some $k$. Thus, we see that embedding
$\V\subset \V_e$ defines an embedding $\slt\subset \g$. Rigidity of
$A$ also implies that the multiplication map $V_{10}\ttt V_{10}\to
V_0$ is non-zero, which implies that the restriction of the commutator
in $\V_e[1]=\g$ to $L_{10}\ttt L_{10}\to \slt$ is non-zero. Now we can
use the following lemma.

\begin{lemma}
Let $\g$ be a finite-dimensional Lie algebra which contains a
subalgebra isomorphic to $\slt$ and as a $\slt$-module,  
$$
\g\simeq \slt\oplus L_{10}. 
$$
If, in addition, restriction of the commutator map $[\, , \, ]\colon
L_{10}\ttt L_{10}\to \slt$ is non-zero, then $\g\simeq G_2$.
\end{lemma}
\begin{proof}
It is easy to see that in such a situation, $\g$ must be simple
(indeed, the only possible ideals are $\slt$ and $L_{10}$, and none of
them is an ideal). But the only 14-dimensional simple Lie algebra is
$G_2$. 
\end{proof}

Therefore, embedding $\V\subset \V_e$ gives rise to an embedding
$\slt\subset G_2$.  Since the Virasoro central charge is the same for
$\V,\V_e$, this embedding extends to a conformal embedding
$(\slthat)_{28}\subset (\ghat)_k$.  But it is well known (see, e.g,
\cite[Chapter 17]{CFT}) that such a conformal embedding uinique,
namely $(\slthat)_{28}\subset (\widehat{G}_2)_1$.

\end{proof}

\begin{remark}
  Note that the proof of \thref{tade} does not rely on
  Itzykson-Cappelli-Zuber classification.
\end{remark}

\begin{remark}
  Explicit analysis shows that for all $\C$-algebras $A$ given by
  \thref{tade}, for any $X,Y\in \Rep A$ one has $X\tta Y\simeq Y\tta
  X$ (not canonically) even though there is no natural way to define
  braiding on $\Rep A$; thus, the Grothendieck ring $K(A)$ is
  commutative. Moreover, this ring coincides with the so-called
  ``graph algebra'' of the Dynkin diagram (see \cite{CFT} for
  discussion of graph algebras). In fact, many of the matrices and
  constants which naturally appear in this theory (such as matrix of
  $F\colon K\to K(A)$) can be calculated using only the Dynkin
  diagram. This was first suggested by Ocneanu \cite{ocneanu} in
  relation with the theory of subfactors; see,
  e.g., \cite{C} for explicit calculations in $E_6$ case.  This
  relation will be discussed in detail elsewhere.
\end{remark}

For future references, we give here some information about $K(A)$ for
Dynkin diagrams of types $D_{even}, E_6$ and $E_8$. This information
can be easily obtained by direct calculation outlined in the proof of
\thref{tade}; checking which of simple $A$-modules lie in $\Rep^0 A$
is trivial: explicit calculation shows that for each of these
algebras, $\th_A=\id$ and thus we can use  \thref{tssrep0}.

\begin{description}
\item[$D_{2m+2}$] This algebra appears when the level $k=4m$; the
  Coxeter number for $D_{2m+2}$ is $l=k+2=4m+2$. The diagram below shows,
  for each of the simple $A$-modules, its structure as an object of
  $\C$. For brevity, we write $i$ instead of $V_i$; thus, $0+(4m)$
  stands for $V_0\oplus V_{4m}$, etc. Filled circles correspond to
  simple objects which lie in $\Rep^0A$; empty circles are simple
  objects in $\Rep A$ which are not in $\Rep^0 A$.  

\fig{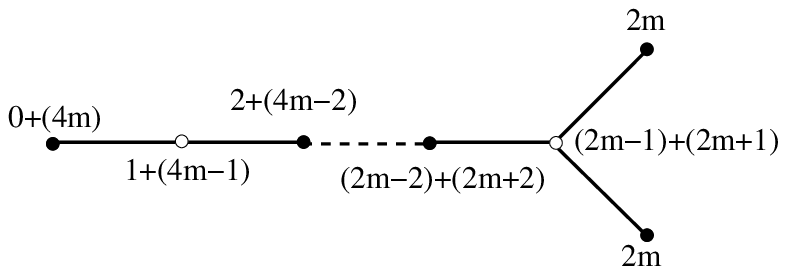}

\item[$E_6$] This algebra appears for $k=10$; the Coxeter number for
  $E_6$ is $l=k+2=12$. All notations are same as before. 

\fig{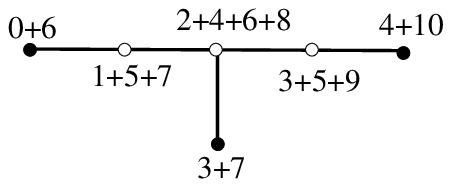}

\item[$E_8$] This algebra appears for $k=28$; the Coxeter number for
  $E_8$ is $l=k+2=30$. 
\medskip

\fig{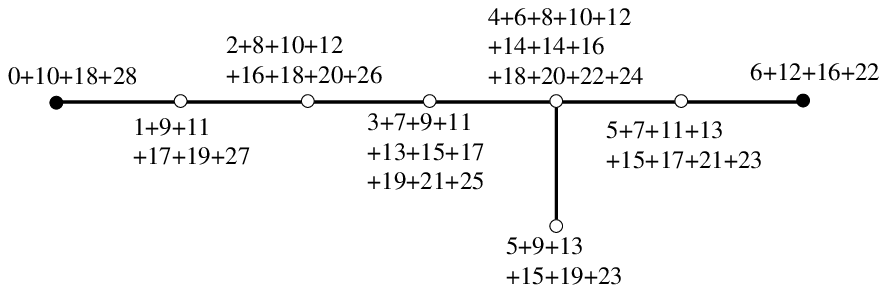}
\end{description}

Note that by \thref{tmodular2}, each of $\C$-algebras $A$ listed in
\thref{tade} gives rise to a modular category $\Rep^0 A$ and thus, a
modular invariant in the sense of conformal field theory. It is easily
checked that these modular invariants coincide with those given by
Cappelli-Itzykson-Zuber classification. Note, however, that our proofs
are completely independent of Cappelli-Itzykson-Zuber classification.

\begin{remark}
  After publication of the first version of this paper, it was
  pointed out  to us that the data given by the figures above had
  previously appeared in the literature in other guises. Most
  importantly, the map $F\colon K\to K(A)$ is a morphism of $K$-modules; in
  particular, this implies that it is an ``intertwiner'' in the sense
  of \cite{FZ}. The explicit formulas for $F$ given above coincide
  with those in Table~1 of \cite{FZ}. However, in the construction in
  \cite{FZ} this map is just one of many possible intertwiners; also,
  they only consider this map at the level of fusion algebras. In our
  approach, $F\colon K\to K(A)$ comes from a functor $F\colon
  \C\to\Rep A$ which is completely determined by the algebra $A$.
\end{remark}

\section{Algebra of type $D_{2n}$}\label{sdtype}

In this section we describe in detail the category of representations
of the algebra $A=\one\oplus\de$ in $\C$, constructed in the previous
section for $k=4m$.

\begin{theorem}\label{tdtype2}

{\ }
\begin{enumerate}
\item 
Simple modules over $A$ are $X_i=V_i \oplus V_{k-i}= A\ttt (V_i),
i=1, \dots, 2m-1$
and two simple modules $X_{2m}^+, X_{2m}^-$, both isomorphic as
objects of $\C$ to $V_{2m}$, with $\mu_{X^+}=\mu_{X^-}\circ p$, where
$p\colon A\to A, p|_{\one}=1, p|_{\de}=-1$. 

\item Tensor product with $F(V_1)=X_1$ is given by 
\begin{equation*}
\begin{aligned}X_1&\tta X_0=X_1,\\
      X_1&\tta X_i\simeq X_{i-1}\oplus X_{i+1},\quad i=1, \dots, 2m-2\\ 
      X_1&\tta X_{2m-1}=X_{2m-2}\oplus X_{2m}^+\oplus X_{2m}^-,\qquad 
      X_1\tta X_{2m}^{\pm}=X_{2m-1}.
    \end{aligned}
  \end{equation*}  
\end{enumerate}
\end{theorem}
Proof is fairly straightforward if we notice that an $A$-module is the
same as an object $V\in \C$ with an isomorphism $\mu\colon \de\ttt
V\isoto V$ such $\mu^2\colon \de\ttt\de\ttt V\to V$ coincides with
$\mu_{\de\de}\ttt \id_V$.

 We also note that formula $F(V)\tta
F(W)\simeq F(V\ttt W)$ defines multiplication in the subring in
$K(A)$ generated by $X_1, \dots, X_{2m-1}$, $(X^+_{2m}+X^-_{2m})$. However,
it does not allow one to determine tensor products involving
$X_{2m}^\pm$. To do so, we need to use the definition.

\begin{theorem}
For $8\mid k$, 

\begin{align*}
X_{2m}^\pm\tta X_{2m}^\pm &\simeq 
X_0\oplus X_4\oplus\dots \oplus X_{2m-4}\oplus X_{2m}^\pm,\\
X_{2m}^\pm\tta X_{2m}^\mp &\simeq 
X_2\oplus X_6\oplus\dots \oplus X_{2m-2}.
\end{align*}
For $k\equiv 4 \mod 8$, 

\begin{align*}
X_{2m}^\pm\tta X_{2m}^\pm &\simeq 
X_2\oplus X_6\oplus\dots \oplus X_{2m-4}\oplus X_{2m}^\mp,\\
X_{2m}^\pm\tta X_{2m}^\mp &\simeq 
X_0\oplus X_4\oplus\dots \oplus X_{2m-2}.
\end{align*}

In particular, $(X^\pm)^*\simeq X^\pm$ for $8\mid k$, and
$(X^\pm)^*\simeq X^\mp$ for $k\equiv 4 \mod 8$,

\end{theorem}

\begin{proof}
By definition, $X\tta Y= (X\ttt Y)/\im (\mu_1-\mu_2)$. As an object of
$\C$, 
$$
X_{2m}^\pm \ttt X_{2m}^\pm=V_{2m}\ttt V_{2m}=V_0\oplus V_2\oplus\dots
\oplus V_k
$$
we need to check which of the modules $V_i$ are in the image of
$\mu_1-\mu_2$. To do so, we use the following lemma.
\begin{lemma}
  Let $n$ be even, $n\le k$ and let $\mu_1, \mu_2\colon \de\ttt
  V_{k-n}\to V_n$ be defined by the compositions
\begin{align*}
\mu_1\colon& \de\ttt V_{k-n}
              \xxto{\id\ttt f} \de\ttt V_{2m}\ttt V_{2m} 
              \xxto{\mu\ttt\id}V_{2m}\ttt V_{2m}
              \xxto{g}V_n\\
\mu_2\colon& \de\ttt V_{k-n}
              \xxto{\id\ttt f} \de\ttt V_{2m}\ttt V_{2m} 
              \xxto{R\ttt \id}V_{2m}\ttt \de\ttt V_{2m}
              \xxto{\id\ttt\mu}V_{2m}\ttt V_{2m}
              \xxto{g}V_n\\
\end{align*}
where $f\colon V_{k-n}\to V_{2m}\ttt V_{2m}$, $g\colon V_{2m}\ttt
V_{2m}\to V_{n}$ and $\mu\colon \de\ttt V_{2m}\to V_{2m}$ are
arbitrary non-zero morphisms. Then $\mu_1=(-1)^{(k-2n)/4}\mu_2$.
\end{lemma}
To prove the lemma, it suffices to consider the identity shown in
\firef{lastfig} and then apply \leref{lfR} to both sides. This proves
the lemma.

\begin{figure}[h]
\begin{equation*}
\fig{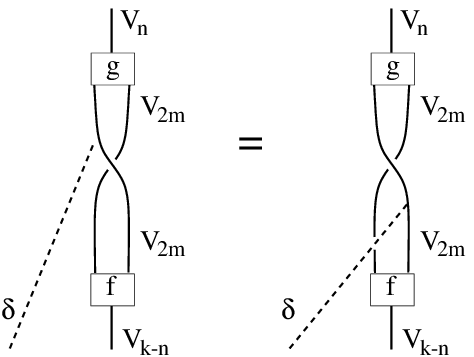}
\end{equation*}
\caption{}\label{lastfig}
\end{figure}

This lemma  implies that for $X^\pm\ttt X^\pm$, $\im
(\mu_1-\mu_2)$ consists of those $V_i$ with $i$ even and $k-2i\equiv 4\mod
8$, while for $X^\pm\ttt X^\mp$, $\im
(\mu_1-\mu_2)$ consists of those $V_i$ with $i$ even and $k-2i\equiv 0\mod
8$. 

This determines the decomposition of $X^\pm\tta X^\pm, X^\pm\tta
X^\mp$ as on object of $\C$. By \thref{tdtype2}, this determines this
tensor product as a representation of $A$ uniquely except for
ambiguity in the choice of the action of $A$ on $V_{2m}$; in other
words, we do not know if $X_{2m}^+$ or $X_{2m}^-$ appears in
decomposition of $X_{2m}^\pm\tta X_{2m}^\pm$. To answer this, note
that we already know enough to deduce that for $8\mid k$,
$(X_{2m}^\pm)^*\simeq X_{2m}^\pm$. Thus, using rigidity we find
$$
\langle X_{2m}^\pm\tta X_{2m}^\pm, X_{2m}^\mp\rangle
=\langle X_{2m}^\pm, X_{2m}^\pm\tta X_{2m}^\mp\rangle=0
$$
since we already know decomposition of $X_{2m}^\pm\tta X_{2m}^\mp$. Similar
arguments show that for $k\equiv 4 \mod 8$, $\langle X^\pm\tta X^\pm,
X^\pm\rangle=0$. This completes the proof of the theorem. 
\end{proof}

\end{document}